\documentclass[11pt]{article}
\title{Volterra-type Ornstein-Uhlenbeck processes in space and time}

\author{
	Viet Son Pham\thanks{Center for Mathematical Sciences, Technical University of Munich, Boltzmannstra\ss e 3, 85748 Garching, Germany, e-mail: vietson.pham@tum.de, carsten.chong@tum.de, URL: www.statistics.ma.tum.de}
	~and
	Carsten Chong$^\ast$
}

\usepackage{latexsym}
\usepackage{amssymb}
\usepackage{amsmath}
\usepackage{amsthm}
\usepackage{mathrsfs}
\usepackage[numbers]{natbib}
\usepackage{stmaryrd}
\usepackage{url}
\usepackage{longtable}
\usepackage{color}
\usepackage{dsfont}
\usepackage[T1]{fontenc}
\usepackage{lmodern}
\usepackage[linktocpage]{hyperref}
\usepackage{graphicx}
\usepackage{float}
\usepackage[latin1]{inputenc}
\usepackage{enumitem}
\setenumerate{leftmargin=*}
\usepackage{amsopn}
\usepackage{comment}

\newcommand{\bfi}{\begin{fig}}
\newcommand{\efi}{\end{fig}}
\newcommand{\btab}{\begin{tab}}
\newcommand{\etab}{\end{tab}}
\newcommand{\barr}{\begin{array}}
\newcommand{\earr}{\end{array}}
\newcommand{\beq}{\begin{equation}}
\newcommand{\eeq}{\end{equation}}
\newcommand{\bdis}{\begin{displaymath}}
\newcommand{\edis}{\end{displaymath}\noindent}

\newcommand{\bbn}{\mathbb{N}}

\newcommand{\bbr}{\mathbb{R}}
\newcommand{\bbc}{\mathbb{C}}
\newcommand{\bbe}{\mathbb{E}}

\newcommand{\bbp}{\mathbb{P}}

\newcommand{\bbf}{\mathbb{F}}
\newcommand{\bbg}{\mathbb{G}}

\newcommand{\bone}{\mathds 1}

\newcommand{\eqd}{\stackrel{\mathrm{d}}{=}}
\newcommand{\halmos}{\quad\hfill $\Box$}
\newcommand{\calh}{{\cal H}}

\newcommand{\calf}{{\cal F}}

\newcommand{\calb}{{\cal B}}

\newcommand{\Omicron}{{\cal O}}

\newcommand{\calm}{{\cal M}}
\newcommand{\calg}{{\cal G}}

\newcommand{\al}{{\alpha}}
\newcommand{\la}{{\lambda}}
\newcommand{\La}{{\Lambda}}
\newcommand{\eps}{{\epsilon}}

\newcommand{\ga}{{\gamma}}
\newcommand{\Ga}{{\Gamma}}
\newcommand{\vp}{{\varphi}}

\newcommand{\si}{{\sigma}}
\newcommand{\om}{{\omega}}
\newcommand{\Om}{{\Omega}}

\newcommand{\var}{{\mathrm{Var}}}
\newcommand{\cov}{{\mathrm{Cov}}}
\newcommand{\corr}{{\mathrm{Corr}}}
\newcommand{\dd}{\mathrm{d}}
\newcommand{\ee}{\mathrm{e}}
\newcommand{\ii}{\mathrm{i}}

\newcommand{\bb}{\mathrm{b}}

\newcommand{\pd}{\partial}

\newcommand{\loc}{\mathrm{loc}}

\newcommand{\DD}{\mathrm{D}}
\newcommand{\CC}{\textcolor{red}}

\newcommand{\Leb}{\mathrm{Leb}}

\newcommand{\pf}{\mathfrak{p}}
\newcommand{\qf}{\mathfrak{q}}

\newcommand{\lpo}{\preceq}
\newcommand{\gpo}{\succeq}

\newcommand{\rpart}{\mathfrak{Re}}

\makeatletter
\newcommand{\opnorm}{\@ifstar\@opnorms\@opnorm}
\newcommand{\@opnorms}[1]{%
  \left|\mkern-1.5mu\left|\mkern-1.5mu\left|
   #1
  \right|\mkern-1.5mu\right|\mkern-1.5mu\right|
}
\newcommand{\@opnorm}[2][]{%
  \mathopen{#1|\mkern-1.5mu#1|\mkern-1.5mu#1|}
  #2
  \mathclose{#1|\mkern-1.5mu#1|\mkern-1.5mu#1|}
}
\makeatother
\newcommand{\acorrf}{\mathrm{acorrf}}	%Autokorrelation
\newcommand\indep{\protect\mathpalette{\protect\independenT}{\perp}}	%independent
\def\independenT#1#2{\mathrel{\rlap{$#1#2$}\mkern2mu{#1#2}}}

\newtheoremstyle{neu}
    {11pt}      % Space above
    {11pt}      % Space below
    {}                  % Body font
    {}          % Indent amount (empty = no indent, \parindent = para indent)
    {\bfseries} % Thm head font
    {}          % Punctuation after thm head
    {1em}  % Space after thm head: " " = normal interword space;
    {\textbf{\thmname{#1}\thmnumber{ #2}\thmnote{ (#3)}}}          % Thm head spec (can be left empty, meaning 'normal')
    
\newtheoremstyle{proof}
    {11pt}      % Space above
    {11pt}      % Space below
    {}                  % Body font
    {}          % Indent amount (empty = no indent, \parindent = para indent)
    {\bfseries} % Thm head font
    {}            % Punctuation after thm head
    {1em}          % Space after thm head: " " = normal interword space;
    {\textbf{\thmname{#1}\thmnote{ #3}.}}          % Thm head spec (can be left empty, meaning 'normal')

\newtheorem{Theorem}{Theorem}[section]
\newtheorem{Corollary}[Theorem]{Corollary}
\newtheorem{Lemma}[Theorem]{Lemma}
\newtheorem{Proposition}[Theorem]{Proposition}
\theoremstyle{neu}
\newtheorem{Definition}[Theorem]{Definition}
\newtheorem{Example}[Theorem]{Example}
\newtheorem{Remark}[Theorem]{Remark}
\newtheorem{Assumption}{Assumption}

\theoremstyle{proof}
\newtheorem{Proof}{Proof}

\newcommand{\bthm}{\begin{Theorem}}
\newcommand{\ethm}{\end{Theorem}}

\newcommand{\bcor}{\begin{Corollary}}
\newcommand{\ecor}{\end{Corollary}}

\newcommand{\blem}{\begin{Lemma}}
\newcommand{\elem}{\end{Lemma}}

\newcommand{\bprop}{\begin{Proposition}}
\newcommand{\eprop}{\end{Proposition}}

\newcommand{\bdf}{\begin{Definition}}
\newcommand{\edf}{\end{Definition}}

\newcommand{\bex}{\begin{Example}}
\newcommand{\eex}{\end{Example}}

\newcommand{\brem}{\begin{Remark}}
\newcommand{\erem}{\end{Remark}}

\newcommand{\bass}{\begin{Assumption}}
\newcommand{\eass}{\end{Assumption}}

\newcommand{\bpr}{\begin{Proof}}
\newcommand{\epr}{\end{Proof}}

\newcommand{\benu}{\begin{enumerate}}
\newcommand{\eenu}{\end{enumerate}}

\newcommand{\bit}{\begin{itemize}}
\newcommand{\eit}{\end{itemize}}

\numberwithin{equation}{section}

\allowdisplaybreaks[3]

\begin{document}

%\date{}

\maketitle

\begin{abstract}
	We propose a novel class of temporo-spatial Ornstein-Uhlenbeck processes as solutions to L\'evy-driven Volterra equations with additive noise and multiplicative drift.
	%We extend the class of L\'evy-driven Ornstein-Uhlenbeck processes to time and space in a novel way. This is achieved by employing techniques from the theory of stochastic Volterra integral equations driven by L\'evy bases. 
	After formulating conditions for the existence and uniqueness of solutions, we derive an explicit solution formula and discuss distributional properties such as stationarity, second-order structure and short versus long memory. Furthermore, we analyze in detail the path properties of the solution process. In particular, we introduce different notions of càdlàg paths in space and time and establish conditions for the existence of versions with these regularity properties. The theoretical results are accompanied by illustrative examples.
\end{abstract}

\vfill

\noindent
\begin{tabbing}
{\em AMS 2010 Subject Classifications:} \= primary: \,\,\,\,\,\, 60G10, 60G17, 60G60, 60H20 \\
\> secondary: \,\,\,60G48, 60G51, 60J75
\end{tabbing}

\vspace{1cm}

\noindent
{\em Keywords:}
ambit process, càdlàg in space and time, Lévy basis, long memory, path properties, second-order structure, space--time modeling, stationary solution, stochastic Volterra equation, Volterra-type Ornstein-Uhlenbeck process

\vspace{0.5cm}

\newpage

\section{Introduction}

In 1908 Langevin \cite{Langevin1908} introduced the following equation as a model for the movement of a particle in a surrounding medium:
\beq\label{OUorig} m \frac{\dd v(t)}{\dd t} = -\la v(t) + \dot W(t).\eeq
Here $v(t)$ denotes the velocity of the considered particle at time $t$, $m$ its mass, $\la$ a friction parameter that accounts for the friction forces acting on the particle, and $\dot W(t)$ the formal derivative of a Brownian motion that governs the random movement of the particle. The solution to the \emph{Langevin equation} \eqref{OUorig} is nowadays called \emph{Ornstein-Uhlenbeck (OU) process}, named after the 1930 paper \cite{Uhlenbeck30}. Employing It\^o's calculus, it is well known that the stochastic differential equation \eqref{OUorig} has a unique solution given by
\beq\label{OUsol}  v(t)= \ee^{-\la t}v(0) + \int_0^t \ee^{-\la(t-s)}\,\dd W(s),\quad t\geq0, \eeq
where $v(0)$ is the initial velocity. Meanwhile, OU processes have found many applications beyond molecular physics, so for example, in stochastic volatility modeling \cite{BN01a}. Moreover, the relatively simple model \eqref{OUsol} has been extended in several directions, for instance, to supOU processes \cite{BN01},  generalized OU processes \cite{Behme_diss} or CARMA processes \cite{Brockwell01}. For all these models, the noise $W$ in \eqref{OUsol} no longer needs to be Gaussian, but can also be a L\'evy process with jumps.

The goal of this article is to extend the class of OU processes to space and time and obtain a temporo-spatial statistical model that is both flexible for modeling purposes and analytically tractable. In fact, extensions of OU processes to several parameters already exist: see \cite{BN04, Brockwell16, Graversen11, Nguyen15} and the first chapter of \cite{Walsh86} for various approaches. They all have in common that they start from \eqref{OUsol} and generalize this formula to several parameters. In this way the main stochastic properties of the one-parameter OU process are preserved because the structure of the exponential kernel is kept. However, in these multi-parameter generalizations, the relation to the original differential equation \eqref{OUorig} is no longer clear. And this is exactly the starting point of our present work: we consider temporo-spatial extensions of OU processes based on the differential equation \eqref{OUorig}. As we shall see, we recover some of the aforementioned extensions as particular cases of our approach. 

The rationale behind our approach is based on two properties of Equation~\eqref{OUorig} which are characteristic to the Langevin equation and which we want to maintain in our generalization. First, the noise $W$ is \emph{additive}, that is, its effect on the process $v$ is independent of the latter. Second, the drift term of $v$ is a scalar \emph{multiple} of its current value. Based on these observations, we propose the following model as a first step towards a temporo-spatial version of \eqref{OUorig}:
\beq\label{TSOU1} \dd X(t,x)=-\la X(t,x)\, \dd t + \int_{\bbr^d} g(x-y)\,\La(\dd t,\dd y),\quad t\geq 0,\ x\in\bbr^d,\eeq
subjected to some initial condition $X(0,x)$ for $x\in\bbr^d$. Here $\La$ is a \emph{homogeneous L\'evy basis} on $[0,\infty)\times \bbr^d$ (also called an infinitely divisible independently scattered random measure or a L\'evy sheet in the literature), which can be thought of as a multi-parameter analogue of a L\'evy process. The function $g$ is such that the integral in \eqref{TSOU1} makes sense, and the differential operation $\dd$ on the left-hand side of \eqref{TSOU1} is taken with respect to time $t$. It is immediate to see that for each fixed $x$, the process $t\mapsto \int_{\bbr^d} g(x-y)\,\La([0,t],\dd y)$ is a L\'evy process and, as a result, the process $t\mapsto X(t,x)$ an OU process in time. Therefore, \eqref{TSOU1} defines a system $(X(\cdot,x))_{x\in\bbr^d}$ of dependent OU processes. Furthermore, if $g$ is strictly positive, every jump of $\La$ induces simultaneous jumps of all OU processes $X(\cdot,x)$, so innovations of $\La$ propagate with an infinite speed through the system. In order to include the case of finite propagation speed, we modify \eqref{TSOU1} by allowing $g$ to depend on time as well:
\beq\label{TSOU2} X(t,x)=X(0,x)-\la \int_0^t X(s,x)\, \dd s + \int_0^t\int_{\bbr^d} g(t-s,x-y)\,\La(\dd s,\dd y),\quad t\geq 0,\ x\in\bbr^d.\eeq
The integral form here is preferable because the kernel is time-dependent. Even more, we can allow the drift term to depend on different neighboring sites, if required, with time delay:
\beq\label{TSOU3} X(t,x)=X(0,x)-\int_0^t \int_{\bbr^d} \mu(t-s,x-y)X(s,y)\, \dd (s,y) + \int_0^t\int_{\bbr^d} g(t-s,x-y)\,\La(\dd s,\dd y),\eeq
where $\mu$ is another kernel function. Unfortunately, model \eqref{TSOU3} no longer contains \eqref{TSOU2} as a special case, but we can remedy this problem by taking a measure $\mu$ instead of a function. So the final class of processes we consider in this article is
\beq\label{TSOU4} X(t,x)=X(0,x)-\int_0^t \int_{\bbr^d} X(t-s,x-y)\, \mu(\dd s,\dd y) + \int_0^t\int_{\bbr^d} g(t-s,x-y)\,\La(\dd s,\dd y).\eeq
If $\mu$ has a density with respect to the Lebesgue measure on $\bbr_+\times\bbr^d$, the last model reduces to \eqref{TSOU3}; if $\mu=-\la \Leb_{\bbr_+}\otimes \delta_{0,\bbr^d}$ (where $\Leb_{\bbr_+}$ is the Lebesgue measure on $\bbr_+$ and $\delta_{0,\bbr^d}$ the Dirac measure in the origin of $\bbr^d$), then it reduces to \eqref{TSOU2}. Henceforth, we refer to \eqref{TSOU4} as the \emph{Volterra-type Ornstein-Uhlenbeck (VOU) equation} and any solution to it as a \emph{VOU process}. We call $\mu$ the \emph{drift measure} and $g$ the \emph{noise propagation function}.

VOU processes have connections to other classes of space--time processes. They are submodels of the general class of \emph{ambit fields} (we refer to the recent survey \cite{BN16}), which have applications in several areas such as turbulence, finance and biological growth modeling. Furthermore, VOU processes are, as we shall show, generalizations of the OU$_\wedge$ process considered in \cite{BN04, Nguyen15}, which has been applied to data of radiation anomalies in the latter reference. In contrast to the OU$_\wedge$ process, the autocorrelation of a VOU process may exhibit long-range dependence and jumps may occur simultaneously for all space locations. In addition to that, VOU processes are solutions to \emph{stochastic Volterra equations} \cite{Chong16, Chong16b} with multiplicative drift and additive noise. The VOU model that we investigate in Example~\ref{Ex2} below is a generalization of the stochastic wave equation in dimension $1$, see \cite{Dalang09}. 

The remaining paper is devoted to a probabilistic analysis of the VOU model and is organized as follows. In Section~\ref{prelim} we recall the necessary background on L\'evy bases and present some results on deterministic Volterra equations that we need throughout the paper. The main theorems of Section~\ref{soltheory}, Theorems~\ref{localsol} and \ref{stationarysol}, give sufficient conditions on $\La$, $g$ and $\mu$ to ensure the existence and uniqueness of solutions to the VOU equation \eqref{TSOU4} and their convergence in distribution as time tends to infinity, respectively. The solution to \eqref{TSOU4} can be expressed as an explicit stochastic convolution integral, revealing the interplay between the theory of deterministic Volterra equations and stochastic convolutions. Section~\ref{soltheory} concludes with a detailed investigation of two Examples~\ref{Ex1} and \ref{Ex2}. In Section~\ref{distprop} we first summarize key distributional properties of the VOU process in Proposition~\ref{gcf} and Corollary~\ref{acf} before we discuss conditions for the VOU model to exhibit temporo-spatial short- or long-range dependence. As Examples~\ref{lrd-ex1} and \ref{lrd-ex2} demonstrate, long memory in \eqref{TSOU4} can arise both through a drift measure with slow decay in time and a measure with slow decay in space. Section~\ref{pathprop} examines the path regularity of the VOU model. If the noise $\La$ is Gaussian, H\"older continuous sample paths can be obtained under mild assumptions, see Theorem~\ref{contversion}. When the noise exhibits jumps, the path properties of the VOU process are basically dictated by the noise propagation function $g$. If it is sufficiently smooth, we show in Theorem~\ref{tempcadlag} that the VOU process has a \emph{t-càdlàg} version, see Definition~\ref{deftempcadlag}. If $g$ is discontinuous, we only have results if the spatial dimension is $1$ and $g=\bone_{-A} h$, where $A$ is a ``triangular'' ambit set in space--time and $h$ is smooth enough. In this case, by Theorem~\ref{pocadlag}, the VOU process has a version which is càdlàg with respect to the triangular shape of $A$, see Definition~\ref{defpocadlag} for a precise statement. Section~\ref{sectproofs} contains the proofs of the main results. Appendix~\ref{AppA} proves Proposition~\ref{thmdetsolution} regarding deterministic Volterra equations, Appendix~\ref{AppB} gives several examples for resolvent measures, and Appendix~\ref{AppC} lists some of their integrability properties.

\section{Preliminaries}\label{prelim}

\subsection{L\'evy bases}\label{notations}

We consider a complete probability space $(\Omega,\calf,\bbp)$ that supports a \emph{homogeneous L\'evy basis} $\La$ on $I\times\bbr^d$ where, depending on the context, $I=\bbr_+=[0,\infty)$ or $I=\bbr$. That is, we assume that $(\La(A))_{A\in\calb_\bb(I\times\bbr^d)}$ is a collection of random variables indexed by bounded Borel subsets of $I\times\bbr^d$ such that for all such $A$ we have
\begin{align}
	\La(A)=&~b \Leb_{I\times\bbr^d}(A) + \sigma W(A) + \int_I\int_{\bbr^d}\int_\bbr \bone_A(t,x)z\bone_{\{|z|\leq 1\}}\,(\pf-\qf)(\dd t,\dd x,\dd z)\nonumber \\
	&~+ \int_I\int_{\bbr^d}\int_\bbr \bone_A(t,x)z\bone_{\{|z|> 1\}}\,\pf(\dd t,\dd x,\dd z),\label{levybasis}
\end{align}
where 
\bit
	\item $b\in\bbr$ and $\si\in\bbr_+$ are constants,
	\item $W$ is Gaussian white noise on $I\times \bbr^d$ such that $W(A)$ has variance $\Leb_{I\times\bbr^d}(A)$ (see e.g. Chapter~I of \cite{Walsh86}),
	\item $\pf$ is a Poisson random measure on $I\times\bbr^d\times \bbr$ with intensity measure $\qf(\dd t,\dd x,\dd z) = \Leb_{I\times\bbr^d}\otimes\nu$, where $\nu$ is a L\'evy measure on $\bbr$ (see e.g. Definition~II.1.20 in \cite{Jacod03}).
\eit
The triplet $(b,\si^2,\nu)$ is referred to as the \emph{characteristics} of $\La$.
Analogously to the L\'evy-It\^o decomposition of L\'evy processes, a L\'evy basis is the sum of a deterministic part, a Gaussian part, a compensated sum of small jumps, and a large jumps part. For more information about the meaning of the integrals with respect to $\pf$ or $\pf-\qf$, we refer to Chapter II of \cite{Jacod03}. If $\int_\bbr |z|\bone_{\{|z| >1 \}}\,\nu(\dd z)<\infty$, we define $b_1:=b+\int_\bbr z\bone_{\{|z| >1 \}}\,\nu(\dd z)$ as the \emph{mean} of $\La$. Similarly, if $\int_\bbr |z|\bone_{\{|z| \leq1 \}}\,\nu(\dd z)<\infty$, we say that $\La$ has \emph{jumps of finite variation} and define the \emph{drift} of $\La$ as $b_0:=b-\int_\bbr z\bone_{\{|z| \leq 1 \}}\,\nu(\dd z)$. Finally, $\La$ is said to be \emph{symmetric} if $b=0$ and $\nu$ is a symmetric measure on $\bbr$.

In this paper we only need Wiener-type stochastic integrals with respect to L\'evy bases since only deterministic integrands will appear. This theory is classic \cite{Rajput89} and we only summarize the most important results we need.
\bprop\label{RR} Suppose that $g\colon I\times\bbr^d\to \bbr$ is a measurable function. The \emph{stochastic integral of $g$ with respect to $\La$}, denoted by
\[ \int_{I\times\bbr^d} g\,\dd \La = \int_{I\times\bbr^d} g(t,x)\,\La(\dd t,\dd x) = \int_I \int_{\bbr^d} g(t,x)\,\La(\dd t,\dd x) \]
either way, is well defined as a limit in probability of approximating simple integrals in the sense of \cite{Rajput89} if and only if 
\begin{enumerate}
	\item $\displaystyle\int_{I\times\bbr^d}\left| b g(t,x)+\int_\bbr (zg(t,x)\bone_{\{ |zg(t,x)|\leq 1 \}}-g(t,x)z\bone_{\{ |z|\leq 1 \} })\,\nu(\dd z)\right|\, \dd (t,x)<\infty$,
	\item $\displaystyle\int_{I\times\bbr^d} \sigma^2|g(t,x)|^2\, \dd (t,x)<\infty$,
	\item $\displaystyle\int_{I\times\bbr^d }\int_\bbr (1\wedge |zg(t,x)|^2)\,\nu(\dd z)\,\dd (t,x)<\infty$.
\end{enumerate}
In this case, the stochastic integral $\int_{I\times\bbr^d} g\,\dd \La$ has an infinitely divisible distribution with characteristic triplet $(b_g,\sigma_g^2, \nu_g)$ given by
%\[ \Phi\left(\int_{S} h\ \dd\La\right)(u)=\exp\left\lbrace iu b_h-\frac{1}{2}u^2C_h+\int_\bbr(e^{iuz}-1-iu\tau(z))\nu_h(\dd z)\right\rbrace , \]
%where
\begin{itemize}
	\item $b_g=\displaystyle\int_{I\times\bbr^d}( b g(t,x)+\int_\bbr (zg(t,x)\bone_{\{ |zg(t,x)|\leq 1 \}}-g(t,x)z\bone_{\{ |z|\leq 1 \} })\,\nu(\dd z))\, \dd (t,x)$,
	\item $\sigma^2_g= \displaystyle\int_{I\times\bbr^d}\sigma^2 |g(t,x)|^2\,\dd (t,x)$,
	\item $\nu_g(B)=\displaystyle \int_{I\times\bbr^d} \int_\bbr \bone_{\{ g(t,x)z\in B \}} \,\nu(\dd z)\,\dd(t,x)$ for any Borel set $B\in\calb(\bbr)$.
\end{itemize}
\eprop
A set of sufficient conditions for the integrability of $g$ with respect to $\La$, which are typically easier to check in practice, is given in Lemma~\ref{integrability-simple}. 

\subsection{Deterministic Volterra equations}\label{VE}

We summarize those results on deterministic convolutional Volterra equations that will be useful in the following sections. The monograph \cite{Gripenberg90} is an excellent reference for single parameter Volterra equations. In our temporo-spatial setting, the equation of interest is given by
\begin{equation}\label{eqnvolterra}
X(t,x)=F(t,x)+\int_0^t\int_{\bbr^d}X(t-s,x-y)\,\mu(\dd s,\dd y),\quad (t,x)\in\bbr_+\times\bbr^d,
\end{equation}
where $X$ is the unknown function, $\mu$ is a signed Borel measure on $\bbr_+\times\bbr^d$ and $F\colon \bbr_+\times\bbr^d\to\bbr$ is a measurable \emph{forcing function}. %, are given\footnote{By $\int_0^t\int_{\bbr^d}X(t-s,x-y)\mu(\dd s,\dd y)$ we mean $\int_{[0,t]\times\bbr^d}X(t-s,x-y)\mu(\dd s,\dd y)$. We choose this notation in order to emphasize the distinction between the time component and the space component.}. In this case a rather complete theory exists, which we present in line with Gripenberg \cite{Gripenberg90}. Note that $t$ is usually interpreted as a time parameter, likewise $x$ its the spatial parameter.\\
Actually, the VOU equation \eqref{TSOU4} is exactly of the form \eqref{eqnvolterra}, except that the forcing function $F$ is stochastic. Therefore, understanding the solution theory to the deterministic problem \eqref{eqnvolterra} is crucial to solving the VOU equation \eqref{TSOU4}.

Before we proceed to Equation~\eqref{eqnvolterra}, let us fix some terminology.
For any Borel subset $S\subseteq\bbr^{d+1}$ we denote by $M(S)$ the space of all signed complete Borel measures on $S$ with finite total variation. As a matter of fact, $M(S)$ becomes a Banach space when equipped with the total variation norm $\|\mu\|:=|\mu|(S)$, where $|\mu|$ is the total variation measure of $\mu$.
We also introduce the notation $M_{\loc}(\bbr_+\times\bbr^d)$ for signed measures on $\bbr_+\times\bbr^d$ which belong to $M([0,T]\times\bbr^d)$ when restricted to $[0,T]\times\bbr^d$ for all positive $T$.
	
Similarly, for $p\in(0,\infty]$, the space $L_{\loc}^p (\bbr_+\times\bbr^d)$ denotes the collection of all functions $\bbr_+\times\bbr^d\to\bbr$ whose restrictions to $[0,T]\times\bbr^d$ belong to $L^p([0,T]\times\bbr^d):=L^p([0,T]\times\bbr^d, \calb([0,T]\times\bbr^d), \Leb_{[0,T]\times\bbr^d})$ for all $T\in\bbr_+$.

Next, for two measures $\mu, \eta\in M(\bbr^{d+1})$ the \emph{convolution} $\mu*\eta$ is the completion of the measure that assigns to each Borel set $B\subseteq\bbr^{d+1}$ the value 
\beq\label{convmeas} (\mu*\eta)(B)=\int_{\bbr^{d+1}}\eta(B-z)\,\mu(\dd z),\eeq
where $B-z=\{s-z\colon s\in B\}$. Since the function $z\mapsto \eta(B-z)$ is Borel measurable and bounded, the integral \eqref{convmeas} is always well-defined. %\\
If $\mu$ and $\eta$ belong to $M(S)$ for some $S\in\calb(\bbr^{d+1})$, we first extend $\mu$ and $\eta$ to $\bbr^{d+1}$ by setting $\bar\mu(B)=\mu(B\cap S)$ and $\bar\eta(B)=\eta(B\cap S)$, then obtain $\bar\mu*\bar\eta$ as above and finally define the convolution $\mu*\eta$ as the restriction of $\bar \mu * \bar\eta$ to $S$. %The convolution of  two measures $\mu, \eta\in M([0,T]\times\bbr^d)$ is defined analogously.
It is customary to write  $\mu^{*0}=\delta_{0,\bbr^{d+1}}$ and $\mu^{*j}=\mu*\mu^{*(j-1)}$ for $j\in\bbn$. % the $(j-1)-$fold convolution of $\mu$ by itself, %if this exists. 
In a similar way, if $\mu\in M(\bbr_+\times\bbr^d)$ and $h\colon\bbr_+\times\bbr^d\to\bbr$ is a measurable function, we define \emph{the convolution of $h$ with respect to $\mu$} as the function $h*\mu=\mu*h$ that is given by
\[ (h*\mu)(t,x)=(\mu*h)(t,x):=\int_0^t \int_{\bbr^d}h(t-s,x-y)\,\mu(\dd s,\dd y),\]
which is defined for those $(t,x)\in\bbr_+\times\bbr^d$ for which the integral exists. Note that the measures considered in this paper may have atoms. Hence we use the convention that integrals over an interval always include the endpoints.

The next theorem is the key result from the theory of convolutional Volterra equations that we need in Section~\ref{soltheory}. It determines conditions under which \eqref{eqnvolterra} has a unique solution. By a \emph{solution} to \eqref{eqnvolterra} we understand a measurable function $X\colon \bbr_+\times\bbr^d \to \bbr$ such that \eqref{eqnvolterra} holds for (Lebesgue-)almost all $(t,x)\in\bbr_+\times\bbr^d$. Two solutions are identified if they agree almost everywhere on $\bbr_+\times\bbr^d$.

\bprop\label{thmdetsolution}
Let $\mu \in M_{\loc}(\bbr_+\times\bbr^d)$ be such that $\mu(\{0\}\times \bbr^d)=0$.
\benu
	\item There exists a unique measure $\rho \in M_{\loc}(\bbr_+\times\bbr^d)$, called the \emph{resolvent of $\mu$}, such that 
	\[ \rho + \mu =\mu*\rho. \]
	\item If $F\in L_{\loc}^p(\bbr_+\times\bbr^d)$ for some $p\in[1,\infty]$, then there exists a unique solution $X\in L_{\loc}^p(\bbr_+\times\bbr^d)$ to \eqref{eqnvolterra}. This solution is given by
		\begin{equation}\label{eqndetsolution}
		X(t,x)=F(t,x)-\int_0^t\int_{\bbr^d} F(t-s,x-y)\,\rho(\dd s,\dd y),\quad (t,x)\in\bbr_+\times\bbr^d,
		\end{equation}
		or in short $X=F-\rho*F$, where $\rho$ is the resolvent of $\mu$.
	\item For every $F\in \mathcal{F}(\mu)$, where
	\[ \mathcal{F}(\mu):=\{F \text{ measurable}\colon  |\mu|*|F|, |\rho|*|F|, (|\mu|*|\rho|)*|F|, (|\mu|*|\rho|^{*2})*|F|<\infty  \text{ a.e.}\}, \]
	%Let $\mathcal{L}(\mu)$ be the space of all measurable functions $F\colon \bbr_+\times\bbr^d\to\bbr$ such that $|\mu|*|F|$, $|\rho|*|F|$ and $(|F|*|\rho|)*|\mu|$ are almost everywhere finite. Then \eqref{eqndetsolution} is the unique solution to \eqref{eqnvolterra} in $\mathcal{L}(\mu)$.
	the function $X$ in \eqref{eqndetsolution} is the unique solution to \eqref{eqnvolterra} in the space 
	\beq\label{Lmu} \mathcal{L}(\mu):=\{L \text{ measurable}\colon  |\mu|*|L|, |\rho|*|L|, (|\mu|*|\rho|)*|L|<\infty  \text{ a.e.}\}.\eeq
	%If $F\colon\bbr_+\times\bbr^d\to\bbr$ is a measurable function such that $F*\rho$, $F*\mu$ and $(F*\rho)*\mu$ are almost everywhere well defined, then \eqref{eqndetsolution} still defines a solution to \eqref{eqnvolterra}. It is unique among those functions with the same stated properties as $F$. 
	\eenu
\eprop
A proof of this theorem, together with some properties and examples of convolutions and resolvents, is given in the Appendix. Note that in a Banach algebra framework the resolvent is also called \emph{quasi-inverse} (cf. Section 2.1 of \cite{Palmer94}).

\section{Solution to the VOU equation}\label{soltheory}

In this section we prove existence and uniqueness of solutions to the VOU equation under general assumptions. In fact, we consider a slightly more general equation than \eqref{TSOU4}, namely
\beq\label{VOU} X(t,x)=V(t,x) + \int_0^t\int_{\bbr^d} X(t-s,x-y)\,\mu(\dd s, \dd y) + \int_0^t\int_{\bbr^d} g(t-s,x-y)\,\La(\dd s,\dd y)\eeq
for $(t,x)\in\bbr_+\times\bbr^d$, where $\mu$ is the drift measure, $g$ is the noise propagation function, $V$ is a measurable stochastic process and $\La$ is a homogeneous L\'evy basis. As usual, we say that a process $\tilde X$ is a version of the process $X$ on $\bbr_+\times\bbr^d$ if for every $(t,x)\in\bbr_+\times\bbr^d$ we have $\tilde X(t,x)=X(t,x)$ almost surely.

%Henceforth, we will refer to $\mu$ as the \emph{drift measure} and to $g$ as the \emph{noise propagation function}. 
%Of course, for \eqref{VOU} to make sense, the stochastic integral of $g$ with respect to $\La$ must be well defined, so $g$ must satisfy the conditions listed in Theorem~\ref{RR}. In general, they cannot be further simplified, but it is possible to formulate sufficient conditions which are easier to verify in practice. 
%Under mild regularity conditions we are able to prove the existence and uniqueness of solutions to the VOU equation \eqref{VOU} . Moreover, this solution can be stated explicitly in terms of the resolvent measure of $\mu$.
\begin{Theorem}\label{localsol} Let $\mu$ be a measure in $M_{\loc}(\bbr_+\times\bbr^d)$ with $\mu(\{0\}\times\bbr^d) =0$ and resolvent $\rho$, %whose realizations almost surely belong to $\mathcal{L}(\mu)$ 
$g\colon\bbr_+\times\bbr^d\to\bbr$ be a measurable function
and $\La$ be a homogeneous L\'evy basis on $\bbr_+\times\bbr^d$ with characteristics $(b,\si^2,\nu)$. We assume that
\beq\label{momLa} \int_\bbr \left(|z|^\alpha \bone_{\{|z|>1\}}+|z|^\beta \bone_{\{|z|\leq 1\}}\right)\,\nu(\dd z)<\infty\eeq
and $g\in L^\al_\loc(\bbr_+\times\bbr^d)\cap L^\beta_\loc(\bbr_+\times\bbr^d)$ for some $\alpha\in(0,1]$, and some $\beta\in[1,2]$ if $\si=0$, and $\beta=2$ if $\si\neq0$. 

If $\al<1$, we further suppose that there exists a submultiplicative weight function $\vp\colon \bbr^d\to\bbr_+$ (that is, a measurable function with $\vp(0)=1$ and
$\vp(x+y)\leq \vp(x)\vp(y)$ for $x,y \in\bbr^d$,
such that $\vp$ is locally bounded and locally bounded away from zero) satisfying $\vp^{-\al}\in L^1(\bbr^d)$, $\vp(|\mu|*|g|)\in L^\infty_\loc(\bbr_+\times\bbr^d)$ and $\vp\mu\in M_\loc(\bbr_+\times\bbr^d)$ (where $(\vp\mu)(\dd t,\dd x):=\vp(x)\,\mu(\dd t,\dd x)$). 
 
Then there exists a measurable version of the process
\[ \int_0^t\int_{\bbr^d} g(t-s,x-y)\,\La(\dd s,\dd y),\quad (t,x)\in\bbr_+\times\bbr^d, \]
in $\calf(\mu)$, and for this version and every measurable process $V$ with almost all paths in $\mathcal{F}(\mu)$, Equation~\eqref{VOU} has a solution with almost all paths in $\mathcal{L}(\mu)$ as in \eqref{Lmu}. A version of this solution is given by 
\beq\label{solformula}
	X(t,x)=V(t,x)-(\rho*V)(t,x)+\int_0^t\int_{\bbr^d} (g-\rho*g)(t-s,x-y)\,\La(\dd s,\dd y),\quad (t,x)\in\bbr_+\times\bbr^d,
\eeq 
or in short $X=V-\rho*V+(g-\rho*g)*\La$.

Moreover, this solution is unique in the sense that for any other solution $\tilde X$ with almost all paths in $\mathcal{L}(\mu)$, we have that almost surely, the paths of $X$ and $\tilde X$ are equal almost everywhere on $\bbr_+\times\bbr^d$.
\end{Theorem}

\begin{Remark}
	\begin{enumerate}
	\item A simple sufficient condition for $V$ to have paths in $\mathcal{F}(\mu)$ is, for example, when the function $(t,x)\mapsto \bbe[|V(t,x)|]$ belongs to $L^\infty_\loc(\bbr_+\times\bbr^d)$. Another would be that the paths of $V$ almost surely belong to $L^p_\loc(\bbr_+\times\bbr^d)$ for some $p\in[1,\infty]$.
	\item Regarding the case $\al\in(0,1)$, typical examples for submultiplicative weight functions include $\vp(x)=(1+|x|)^\eta(\log(\ee \vee |x|))^\ga$ and $\vp(x)=\exp(|x|^\ga)$ for $\eta,\ga\geq0$. The reason why we impose additional conditions when $\al$ is smaller than one is that in our proof we have to ensure $(|\rho|*|g|)^\alpha\in L^1_\loc(\bbr_+\times\bbr^d)$. Instead of formulating conditions on $\rho$, which may not be known explicitly, the assumptions in Theorem~\ref{localsol} are solely on $\mu$.
	\end{enumerate}
	\halmos
\end{Remark}

We also remark that the condition in \eqref{momLa} on the L\'evy measure $\nu$ are not necessary in general. For instance, it is well known that the OU process is defined for all L\'evy processes without any restrictions on $\nu$. But this is different to our case because the spatial coordinate is in the non-compact space $\bbr^d$ and the noise propagation function $g$ is not necessarily bounded. For given classes of $g$ and $\mu$ it may be possible to relax the assumptions of Theorem~\ref{localsol}. But given that these are already general enough to cover most practical cases, we refrain from doing so.

Having clarified the local existence of solutions to \eqref{VOU}, our next aim is to investigate their long-term behavior and the existence of stationary solutions.
A stochastic process $X$ on $\bbr_+\times\bbr^d$ is called \emph{strictly stationary} if for every $n\in\bbn$ and $(\tau,\xi),(t_1,x_1),\ldots,(t_n,x_n)\in\bbr_+\times\bbr^d$ the distributions of $(X(t_1,x_1),...,X(t_n,x_n))$ and $(X(t_1+\tau,x_1+\xi),...,X(t_n+\tau, x_n+\xi))$ are equal.

\bthm\label{stationarysol} Let $\La$ be a homogeneous L\'evy basis on $\bbr\times\bbr^d$, the conditions of Theorem~\ref{localsol} be valid with $\al,\beta\in(0,2]$ and additionally $g \in L_{\loc}^{1}(\bbr_+\times\bbr^d)$ such that the process $X$ as given in \eqref{solformula} is the solution to \eqref{VOU}. Moreover, we assume the following hypotheses: 
\benu
	\item For all $x_1,\ldots,x_n\in\bbr^d$ we have that
	\begin{equation*}			
	%\label{convV} 
	\big(V(t,x_1)-(\rho*V)(t,x_1),\ldots, V(t,x_n)-(\rho*V)(t,x_n)\big) \stackrel{\dd}{\longrightarrow} \big(F_\infty(x_1),\ldots,F_\infty(x_n)\big)
	\end{equation*}	
	as $t\to\infty$
	for some deterministic measurable function $F_\infty\colon\bbr^d\to\bbr$.
	\item We have $g-\rho*g \in L^{\al}(\bbr_+\times\bbr^d)\cap L^\beta(\bbr_+\times\bbr^d)$.
	\item Conditions \eqref{infassym} and \eqref{0assym} are satisfied with $K(A):=A^{1-\alpha}$ and $k(a):=a^{1-\beta}$.
\eenu

Then we have for all $n\in\bbn$ and $x_1,\ldots,x_n\in\bbr^d$ that
\begin{equation*}
	%\label{convX}
	\big(X(t,x_1),\ldots,X(t,x_n)\big) \stackrel{\mathrm{d}}{\longrightarrow} \big(X_\infty(x_1),\ldots, X_\infty(x_n)\big),\quad t\to\infty,	
\end{equation*}	
where $X_\infty$ is the spatial process
\begin{equation*}
	%\label{Xinfty} 
	X_\infty(x):=F_\infty(x)+\int_0^\infty \int_{\bbr^d} (g-\rho*g)(s,x-y)\, \La(\dd s,\dd y),\quad x\in\bbr^d.	
\end{equation*}

Furthermore, if $g$ is integrable with respect to $\La$ and $V$ is independent of $\La$ with the same finite-dimensional distributions as
\beq\label{Vdist} \int_{-\infty}^0\int_{\bbr^d} g(t-s,x-y)\,\La(\dd s,\dd y),\quad (t,x)\in\bbr_+\times\bbr^d, \eeq
then $X$ is a strictly stationary process on $\bbr_+\times\bbr^d$. In particular, if $V$ equals the process in \eqref{Vdist}, $X$ can be written as the two-sided strictly stationary process
\beq\label{twosidedX} X(t,x)=\int_{-\infty}^t\int_{\bbr^d} (g-\rho*g)(t-s,x-y)\,\La(\dd s,\dd y),\quad(t,x)\in\bbr\times\bbr^d.\eeq
\ethm

We give various possibilities of how to ensure the requirements in Theorem~\ref{stationarysol}.
\blem\label{statcondLa} In the following cases, condition \eqref{infassym} with $K(A):=A^{1-\alpha}$ (resp. \eqref{0assym} with $k(a):=a^{1-\beta}$) is already implied by \eqref{momLa}:
\benu
\item $\al \in (0,1]$ (resp. $\beta\in[1,2]$).
\item $\al\in(1,2]$ and $b_1=0$ (resp. $\beta\in(0,1)$ and $b_0=0$).
\item $\La$ is symmetric.
\eenu
\elem

\blem\label{rhogint}  For $p\in[1,\infty]$ we have $g-\rho*g\in L^p(\bbr_+\times\bbr^d)$ under each of the following assumptions:
\benu
	\item $\rho\in M(\bbr_+\times\bbr^d)$ and $g\in L^p(\bbr_+\times\bbr^d)$.
	\item $\rho\in M(\bbr_+\times\bbr^d)$ satisfies $\rho(\bbr_+\times\bbr^d)=1$ and there exists a constant $g_\infty\in\bbr$ such that $g-g_\infty \in L^p(\bbr_+\times\bbr^d)$.
	\item For some $q,s\in [1,\infty]$ satisfying $s^{-1}+q^{-1} = 1+ p^{-1}$ we have $g\in L^p(\bbr_+\times\bbr^d) \cap L^s(\bbr_+\times\bbr^d)$ and $\rho(\dd t,\dd x) = r(t,x)\,\dd (t,x)$ with $r\in L^q(\bbr_+\times\bbr^d)$.
\eenu
\elem

In contrast to Theorem~\ref{localsol}, the conditions imposed in Theorem~\ref{stationarysol} (and also in Lemma~\ref{rhogint}) explicitly depend on the behavior of the resolvent measure $\rho$, instead merely on $\mu$. In fact, there are no general necessary and sufficient conditions for a measure $\mu$ to have a resolvent with certain integrability properties. In Appendix~\ref{AppC} we present several results in this respect.

We conclude this section by investigating two variants of a VOU process sharing the same drift measure as in the classical Ornstein-Uhlenbeck process, namely
\beq\label{OUmu} \mu=-\la\Leb_{\bbr_+}\otimes\delta_{0,\bbr^d}, \eeq
and different choices for the noise propagation function $g$.

\bex[VOU process with infinite speed propagation of noise]\label{Ex1}~

\noindent As a first example we investigate the equation
\begin{equation}\label{eqmodel1}
X(t,x)=-\la\int_0^t X(s,x)\, \dd s + \int_0^t\int_{\bbr^d}\ee^{-\la'|x-y|}\,\La(\dd s,\dd y),\quad (t,x)\in\bbr_+\times\bbr^d,
\end{equation}
where $\la\in\bbr$ and $\la'>0$ and $|\cdot|$ denotes the Euclidean norm in $\bbr^d$. A closer inspection reveals two characteristic features of this model: first, the parameter $\la$ leads to a mean-reverting behavior in time like in the classical OU case if $\la>0$; second, since the noise propagation function $g(t,x)=\ee^{-\la'|x|}$ does not depend on $t$ and is strictly positive, each innovation of $\La$ (and jump if the L\'evy measure is not zero) affects $X(\cdot,x)$ for all $x$ simultaneously. However, as controlled by $\la'$, the impact of an innovation decreases exponentially in the distance between the current location $x$ and the point of origin $y$. For fixed $x$, we further observe that the second summand on the right-hand side of \eqref{eqmodel1} is a L\'evy process, so the solution of \eqref{eqmodel1} is in fact a system $(X(\cdot,x)\colon x\in\bbr^d)$ of dependent classical OU processes. 

Since the resolvent measure of \eqref{OUmu} is $\rho(\dd t,\dd x)=\la \ee^{-\la t}\,\dd t\,\delta_{0,\bbr^d}(\dd x)$ (see Lemma~\ref{resLemma}), a simple calculation yields $(g-\rho*g)(t,x)=\ee^{-\la t -\la'|x|}$. Therefore, as soon as the L\'evy measure $\nu$ of $\La$ satisfies 
\beq\label{assmom} \int_\bbr |z|^\al\bone_{\{ |z|>1 \}}\,\nu(\dd z)<\infty\eeq
for some $\al>0$, we derive from Theorem~\ref{localsol} that the unique solution to \eqref{eqmodel1} is given by
\begin{equation}\label{eqsolmodel1}
X(t,x)=\int_0^t\int_{\bbr^d}\ee^{-\la(t-s)-\la'|x-y|}\,\La(\dd s,\dd y),\quad (t,x)\in\bbr_+\times\bbr^d.
\end{equation}	
Indeed, we can choose $\al$ as above and $\beta=2$ in \eqref{momLa} because $g\in L^p_\loc(\bbr_+\times\bbr^d)$ for all $p\in(0,\infty]$. If $\al<1$, we can take $\vp(x):=(1+|x|)^{(d+1)/\al}$.

Against the background that $X(\cdot,x)$ is an OU process for fixed $x$, it is not surprising that also Theorem~\ref{stationarysol} applies if $\la>0$ (then $g-\rho*g\in L^p(\bbr_+\times\bbr^d)$ for all $p\in(0,\infty]$). In this case the strictly stationary process \eqref{twosidedX} is given by
\beq\label{twosidedX1} X(t,x) = \int_{-\infty}^t \int_{\bbr^d} \ee^{-\la (t-s) + \la'|x-y|}\,\La(\dd s,\dd y),\quad (t,x)\in\bbr\times\bbr^d. \eeq
\halmos\eex

In the previous example innovations of $\La$ at a given site have an instantaneous effect on all other sites. In contrast to this, the next model incorporates a traveling waves mechanism such that a certain amount of time is needed for the propagation of innovations from one to another point in space.

\bex[VOU process with finite speed propagation of noise]\label{Ex2}~

\noindent We consider 
\begin{equation}\label{eqmodel2}
X(t,x)=-\la \int_0^t X(s,x)\, \dd s + \int_0^t\int_{\bbr^d}\bone_{\{|x-y|\leq c(t-s)\}}\ee^{-\la'|x-y|}\,\La(\dd s,\dd y),\quad (t,x)\in\bbr_+\times\bbr^d,
\end{equation}
with parameters $c>0$ and $\la,\la'\in\bbr$. As a result, the time until an innovation of $\La$ at a site $y$ arrives at another site $x$ amounts to $|x-y|/c$. With $g(t,x)=\bone_{\{|x|\leq ct\}}\ee^{-\la'|x|}$, an elementary computation shows that 
\begin{equation*}
	\label{kernel2} (g-\rho*g)(t,x) =\bone_{\{|x|\leq ct\}}\ee^{-\la t -(\la'-\la/c)|x|},\quad(t,x)\in\bbr_+\times\bbr^d. 
\end{equation*}	
Consequently, under assumption \eqref{assmom}, Equation~\eqref{eqmodel2} has the unique solution
\begin{equation}\label{eqsolmodel2}
	X(t,x)=\int_0^t\int_{\bbr^d}\bone_{\{|x-y|\leq c(t-s)\}}\ee^{-\la (t-s)-(\la'-\la/c)|x-y|}\,\La(\dd s,\dd y),\quad (t,x)\in\bbr_+\times\bbr^d.
\end{equation}	

In order to determine the long-term behavior of \eqref{eqsolmodel2}, we can use Fubini's theorem to verify that only for $\la>0$ and $\la'>0$ the integral  
\begin{align*} \int_0^\infty \int_{\bbr^d} (g-\rho*g)^p(t,x)\,\dd(t,x) =&~\frac{2\pi^{d/2}}{\Ga(\frac{d}{2})} \int_0^\infty \ee^{-\la p t} \int_0^{ct} \ee^{-(\la' -\la/c)pr}r^{d-1}\,\dd r\,\dd t\\
=&~\frac{2\pi^{d/2}}{\Ga(\frac{d}{2})} \int_0^\infty \frac{\ee^{-r\la p/c}}{\la p} \ee^{-(\la' -\la/c)pr}r^{d-1}\,\dd r\\
=&~ \frac{2\pi^{d/2}}{\Ga(\frac{d}{2})\la p} \int_0^\infty \ee^{-\la'pr}r^{d-1}\,\dd r
\end{align*}
is finite for $p\in(0,\infty)$. So only in this case, the finite-dimensional distributions of \eqref{eqsolmodel2} converge to that of the process
\beq\label{twosidedX2} X(t,x)= \int_{-\infty}^t\int_{\bbr^d}\bone_{\{|x-y|\leq c(t-s)\}}\ee^{-\la (t-s)-(\la'-\la/c)|x-y|}\,\La(\dd s,\dd y),\quad (t,x)\in\bbr\times\bbr^d.\eeq
We notice that if $\la'=\la/c$, the process $X$ in \eqref{eqsolmodel2} is exactly the so-called OU$_\wedge$ model investigated in \cite{BN04, Nguyen15}. 
\halmos\eex

\begin{Remark}
	At the end of this section we want to highlight a connection to stochastic partial differential equations as studied in \cite{Peszat07}. For this purpose assume in addition to the conditions of Theorem \ref{localsol} that $\La$ has mean zero and a finite second moment. Let $U$ be a Hilbert space such that the embedding of $\calh:=L^2(\bbr^d)$ into $U$ is dense and Hilbert-Schmidt (see Example 14.25 of \cite{Peszat07} for an example of $U$). Then the process $W\colon [0,\infty)\times\calh\to L^2(\Om)$ defined by
	\begin{equation*}
	W(t,\phi):=\int_0^t\int_{\bbr^d}\phi(y)\,\La(\dd s,\dd y)
	\end{equation*}
	is the sum of  a cylindrical Wiener process and an impulsive cylindrical process on $\calh$ (cf. Definitions~7.11 and 7.23 of \cite{Peszat07}). Combining Theorems 7.13 and 7.22 of \cite{Peszat07} we obtain a $U$-valued square-integrable Lévy martingale $L$ satisfying
	\begin{equation*}
	L(t)=\sum_{n\in\bbn}W(t,e_n)e_n
	\end{equation*}
	for any fixed orthonormal basis $(e_n)_{n\in\bbn}$ in $\calh$ (cf. Remark 7.14 of \cite{Peszat07}). 
	
	If we extend the convolution operators
	\begin{equation*}
	S(t)\colon\calh\to\calh,\quad\phi(x)\mapsto\int_{\bbr^d}g(t,x-y)\phi(y)\,\dd y
	\end{equation*}
	onto $U$, we may rewrite
	\begin{align*}
	\int_0^t\int_{\bbr^d} g(t-s,x-y)\,\La(\dd s,\dd y)&=\int_0^t\int_{\bbr^d}\sum_{n\in\bbn}\int_{\bbr^d} g(t-s,x-z)e_n(z)\,\dd z \, e_n(y)\,\La(\dd s,\dd y)\\
	&=\sum_{n\in\bbn}\int_0^t\int_{\bbr^d} g(t-s,x-z)e_n(z)\,\dd z \, \int_{\bbr^d}e_n(y)\,\La(\dd s,\dd y)\\
	&=\sum_{n\in\bbn}\int_0^t S(t-s)e_n \, \dd W(s,e_n)=\int_0^t S(t-s)\,\dd L(s).
	\end{align*}
	From the solution formula \eqref{solformula} (with $V\equiv0$) we see that $X(t,x)$ belongs to $L^2(\bbr^d,(1+|x|^r)^{-1})$ for fixed $t$ and $r>d/2$, where a function $\phi$ is an element of $L^2(\bbr^d,\eta)$ if and only if $\eta\phi\in L^2(\bbr^d)$.
	
	Further assuming that $\mu(\dd t,\dd x)=\nu_t(\dd x)\,\dd t$ for a transition kernel $\nu$ (as in Examples \ref{Ex1} and \ref{Ex2}) and that the mapping
	\begin{equation*}
		Q(t)\colon g(x)\mapsto\int_{\bbr^d}g(x-y)\,\nu_t(\dd y)
	\end{equation*}
	is a linear convolution operator from $L^2(\bbr^d,(1+|x|^r)^{-1})$ into itself, we obtain
	\begin{equation}
		\int_0^t\int_{\bbr^d} X(t-s,x-y)\,\mu(\dd s, \dd y)=\int_0^t\int_{\bbr^d} X(s,x-y)\,\nu_{t-s}(\dd y)\,\dd s=\int_0^t Q(t-s)X(s)\,\dd s,
	\end{equation}
	where $X(t):=X(t,\cdot)$. In short, under the conditions above, the VOU equation
	\begin{equation*}
		X(t,x)=X(0,x)+\int_0^t \int_{\bbr^d} X(t-s,x-y)\, \mu(\dd s,\dd y) + \int_0^t\int_{\bbr^d} g(t-s,x-y)\,\La(\dd s,\dd y)
	\end{equation*}	
	is equivalent to the infinite-dimensional equation
	\begin{equation*}
		X(t)=X(0)+\int_0^t Q(t-s)X(s)\,\dd s+\int_0^t S(t-s)\,\dd L(s).
	\end{equation*}
	In the literature for stochastic partial differential equations several criteria are known for $t\mapsto X(t)$, viewed as a process with values in a Hilbert space, to have continuous or càdlàg sample paths (see for instance Theorem~11.8 of \cite{Peszat07} or Theorem~4.5 and Remark~4.6 of \cite{Tappe12}). By contrast, the random field approach to the VOU equation allows for a detailed analysis of the temporo-spatial path properties of $(t,x)\mapsto X(t,x)$ as in Section~\ref{pathprop}.
	\halmos
\end{Remark}

\section{Distributional properties}\label{distprop}

A convenient tool for characterizing the distribution of temporo-spatial processes is the \emph{generalized cumulant functional} introduced in \cite{BN15}. For the solution process $X$ in \eqref{solformula}, we obtain the following result.
\bprop\label{gcf} Suppose that the conditions of Theorem~\ref{localsol} are satisfied with $V\equiv0$ and that $X$ is the solution to \eqref{VOU} given in \eqref{solformula}. If $m\in M(\bbr_+\times\bbr^d)$ is supported on a compact subset of $\bbr_+\times \bbr^d$, the integral $m[X] := \int_{\bbr_+}\int_{\bbr^d} X(t,x)\,m(\dd t,\dd x)$ is well defined and the generalized cumulant functional of $X$ with respect to $m$ is given by
\begin{align*} 
\log \bbe\left[\ee^{\ii u m[X]} \right] = \ii b_G u - \frac{1}{2}\si^2_G u^2 + \int_\bbr \left(\ee^{\ii u z} - 1 - \ii u z \bone_{\{|z|\leq 1 \}}\right)\, \nu_G(\dd z),\quad u\in\bbr,
\end{align*}
where $(b_G, \si^2_G, \nu_G)$ are the characteristics as given in Proposition~\ref{RR} for the function 
\[ G(s,y) = \int_s^\infty \int_{\bbr^d} (g-\rho*g)(t-s,x-y)\,m(\dd t,\dd x),\quad (s,y)\in\bbr_+\times\bbr^d. \]
\eprop

For example, if one takes $m(\dd t,\dd x)=\theta_1\delta_{(t_1,x_1)} + \dots + \theta_n\delta_{(t_n,x_n)}$, one obtains the cumulant-generating function of $(X(t_1,x_1),\ldots,X(t_n,x_n))$. Based on this, it is also possible to derive the second-order structure for $X$, see Proposition~2 in \cite{BN15} for a proof.
\bcor\label{acf} Suppose that the assumptions of Theorem~\ref{localsol} are satisfied with $V\equiv0$ and that $X$ is the solution process \eqref{solformula}.
\benu
	\item If the assumptions of Theorem~\ref{localsol} hold with $\al=1$, $X(t,x)$ has a finite first moment for all $(t,x)\in\bbr_+\times\bbr^d$ given by  
	\beq\label{meanX} \bbe[X(t,x)]=b_1 \int_0^t\int_{\bbr^d} (g-\rho*g)(s,y)\,\dd(s,y). \eeq
	\item If the assumptions of Theorem~\ref{localsol} are additionally satisfied with $\al=2$, $X(t,x)$ has a finite second moment for all $(t,x)\in\bbr_+\times\bbr^d$ and  
	\beq\label{varX} \var[X(t,x)] = \left(\si^2 + \int_\bbr z^2\,\nu(\dd z)\right) \int_0^t\int_{\bbr^d} (g-\rho*g)^2(s,y)\,\dd (s,y). \eeq 
	Moreover, for $(t,x),(\tau,\xi)\in\bbr_+\times\bbr^d$ we have
	\begin{align} &~\cov[X(t,x),X(t+\tau,x+\xi)] \nonumber\\
	&~\qquad = \left(\si^2 + \int_\bbr z^2\,\nu(\dd z)\right)\int_0^t \int_{\bbr^d} (g-\rho*g)(s,y)(g-\rho*g)(s+\tau,y+\xi)\,\dd(s,y). \label{covX} \end{align}
\eenu
Furthermore, in the setting of Theorem~\ref{stationarysol}, if $X$ is the strictly stationary process \eqref{twosidedX}, then the formulae \eqref{meanX}, \eqref{varX} and \eqref{covX} remain valid if we replace $t$ by $\infty$ on the right-hand sides. 
\ecor

For illustration we calculate the autocorrelation functions for the models in Examples~\ref{Ex1} and \ref{Ex2}.

\begin{Example}[Second-order structure for Example~\ref{Ex1}] ~~

\noindent Under the moment assumptions of Corollary~\ref{acf}, the mean and the autocovariance function of the stationary process $X$ in \eqref{twosidedX1} are given by
\begin{equation}
	\label{meanX1} 
	\bbe[X(t,x)]=b_1\int_{\bbr_+}\int_{\bbr^d}\ee^{-\la s-\la'|y|}\, \dd (s,y) = \frac{2b_1\pi^{d/2}\Ga(d)}{\la(\la')^d\Ga(\frac{d}{2})},\quad (t,x)\in\bbr\times\bbr^d,	
\end{equation}
	and, denoting $m_2:=\si^2 + \int_\bbr z^2\,\nu(\dd z)$,
	\begin{align*}
\cov[X(t,x),X(t+\tau,x+\xi)]=&~m_2 \int_{\bbr_+}\int_{\bbr^d}\ee^{-\la s-\la'|y|}\ee^{-\la(s+\tau)-\la'|y+\xi|}\, \dd(s,y) \nonumber\\
=&~\frac{m_2\ee^{-\la\tau}}{2\la}\int_{\bbr^d}\ee^{-\la'|y|-\la'|y+\xi|}\, \dd y=:\frac{m_2\ee^{-\la\tau}}{2\la}E(\xi). 
	\end{align*}
	The integral $E(\xi)$ is the convolution of the function $f(\xi)=\ee^{-\la' |\xi|}$ with itself in $\bbr^d$. Since the Fourier transform of $f$ is known (see Theorem~I.1.14 in \cite{Stein71}), $E(\xi)$ is the inverse Fourier transform of the function 
	\[ \frac{c_d^2a^2}{(a^2+|x|^2)^{d+1}}, \]
	where $a:=\la'/(2\pi)$ and $c_d:=\Gamma((d+1)/2)\pi^{-(d+1)/2}$.
	Hence, using Theorem~IV.3.3 of \cite{Stein71} and denoting by $J_\al$ and $K_\al$ the Bessel functions of the first kind and the modified Bessel functions of the second kind, respectively, we obtain
	\begin{align*}
	E(\xi) = 2\pi |\xi|^{1-d/2} c_d^2 a^2 \int_0^\infty \frac{J_{d/2-1}(2\pi|\xi|r)r^{d/2}}{(a^2+|r|^2)^{d+1}}\,\dd r = \frac{2\Gamma(\textstyle\frac{d+1}{2})^2}{\Gamma(d+1)}\left(\frac{\lambda'}{2\pi} \right)^{1-d/2}|\xi|^{1+d/2} K_{1+d/2}(\la' |\xi|).
	\end{align*}
	This yields for $(t,x)\in\bbr\times\bbr^d$ and $(\tau,\xi)\in\bbr_+\times\bbr^d$
	\begin{equation*}
		%\label{covX1} 
		\cov[X(t,x),X(t+\tau,x+\xi)] = \frac{m_2\Gamma(\textstyle\frac{d+1}{2})^2}{\la\Gamma(d+1)}\left(\frac{\lambda'}{2\pi} \right)^{1-d/2}\ee^{-\la\tau}|\xi|^{1+d/2} K_{1+d/2}(\la' |\xi|). 	
	\end{equation*}
	Since $\lim_{x\downarrow 0} x^{\al}K_{\al}(x) = 2^{\al-1}\Gamma(\al)$ for $\al\geq0$, the autocorrelation function reads as
	\begin{equation*}
		%\label{acorrf1} 
		\corr[X(t,x),X(t+\tau,x+\xi)] = \frac{(\la')^{1+d/2}}{2^{d/2}\Ga(1+\textstyle\frac{d}{2})}\ee^{-\la\tau}|\xi|^{1+d/2} K_{1+d/2}(\la' |\xi|).	
	\end{equation*} 
	If $d=1$, $d=2$ or $d=3$, this formula reduces to
	\begin{align*}
	\corr[X(t,x),X(t+\tau,x+\xi)] =&~ \ee^{-\la\tau}(\la'|\xi|+1)\ee^{-\la'|\xi|}, &(\tau,\xi)&\in\bbr_+\times\bbr,\\
	\corr[X(t,x),X(t+\tau,x+\xi)] =&~ \textstyle\frac{(\la')^2}{2}\ee^{-\la\tau}|\xi|^2K_2(\la'|\xi|),  &(\tau,\xi)&\in\bbr_+\times\bbr^2,\\
	\corr[X(t,x),X(t+\tau,x+\xi)] =&~ \ee^{-\la\tau}\left(\textstyle\frac{(\la')^2}{3}|\xi|^2 +\la'|\xi|+1\right)\ee^{-\la'|\xi|}, &(\tau,\xi)&\in\bbr_+\times\bbr^3,
	\end{align*}
	%\begin{align*}
	%\cov[X(t,x),X(t+\tau,x+\xi)] =&~ \frac{m_2}{2\la}\ee^{-\la\tau}\frac{\la'|\xi|+1}{\la'}\ee^{-\la'|\xi|}, &(\tau,\xi)&\in\bbr_+\times\bbr,\\
	%\cov[X(t,x),X(t+\tau,x+\xi)] =&~ \frac{m_2\pi}{8\la}\ee^{-\la\tau}|\xi|^2K_2(\la'|\xi|),  &(\tau,\xi)&\in\bbr_+\times\bbr^2,\\
	%\cov[X(t,x),X(t+\tau,x+\xi)] =&~ \frac{m_2\pi}{6\la}\ee^{-\la\tau}\frac{(\la')^2|\xi|^2 +3\la'|\xi|+3}{(\la')^3}\ee^{-\la'|\xi|}, &(\tau,\xi)&\in\bbr_+\times\bbr^3,
	%\end{align*}
	respectively.	
	\halmos
\end{Example}
\begin{Example}[Second-order structure for Example~\ref{Ex2}] ~~
	
	\noindent 
	We obtain the same value as in \eqref{meanX1} for the expectation of \eqref{twosidedX2}: 
	\begin{equation*}
		%\label{meanX2} 
		\bbe[X(t,x)]=b_1\int_{\bbr_+}\int_{\bbr^d}\bone_{\{ |y|\leq cs \}}\ee^{-\la s-(\la' -\la/c)|y|}\, \dd (s,y)=\frac{2b_1\pi^{d/2}\Ga(d)}{\la(\la')^d\Ga(\frac{d}{2})},\quad  (t,x)\in\bbr\times\bbr^d. 	
	\end{equation*}
	Regarding the autocovariance function, a straightforward calculation for $d=1$ shows that
	 \begin{align*}
		&\cov[X(t,x),X(t+\tau,x+\xi)]\\
		&~\qquad = \frac{m_2}{4\la} \ee^{-\la\tau} \ee^{(\la/c-\la')|\xi|} \bigg(  \frac{c}{\la}\left(\ee^{-\la(|\xi|/c-\tau)_+}+\ee^{-2\la(|\xi|/c-\tau)}\ee^{\la(|\xi|/c-\tau)_+}-\ee^{-2\la(|\xi|/c-\tau)}\right)\\
		&~\qquad \quad +\frac{1}{\la'}\ee^{-2\la(|\xi|/c-\tau)_+} \bigg) + \frac{m_2}{4\la} \ee^{-\la\tau}  \ee^{-(\la/c+\la')|\xi|} \left(\frac{1}{\la'}-\frac{c}{\la}\right)					
	 \end{align*}
	 for all $(t,x)\in\bbr\times\bbr$ and $(\tau,\xi)\in\bbr_+\times\bbr$.
	 \halmos
\end{Example}

The autocovariance function in the last example exhibits an exponential decay in both $t$ and $x$, so the corresponding process $X$ has a \emph{short-range dependence} structure. However, as we shall see, under suitable choices of $g$ and $\mu$, it may happen that the autocovariance function is not integrable, i.e.
\begin{equation}\label{longmem}
	\int_0^\infty\int_{\bbr^d} |\cov[X(t,x),X(t+\tau,x+\xi)]|\,\dd(\tau,\xi) = \infty,
\end{equation}	
hence giving rise to models with \emph{long-range dependence}. A first result concerning short- versus long-range dependence is the following.

\bprop\label{SRD} Let $\rho\in M_\loc(\bbr_+\times\bbr^d)$ be the resolvent measure associated to some measure $\mu\in M_\loc(\bbr_+\times\bbr^d)$ with $\mu(\{0\}\times\bbr^d)=0$ and suppose that the L\'evy basis $\La$ has a finite second moment, that is $\int_{\{|z|>1\}} |z|^2\,\nu(\dd z)<\infty$.
\benu
	\item If in addition $g-\rho*g\in L^1(\bbr_+\times\bbr^d)\cap L^2(\bbr_+\times\bbr^d)$, then the process $X$ in \eqref{twosidedX} is well defined, has a finite second moment and 
\[ \int_0^\infty\int_{\bbr^d} |\cov[X(t,x),X(t+\tau,x+\xi)]|\,\dd(\tau,\xi) < \infty. \]
	\item If $\La$ has zero mean, $g-\rho*g \in L^1_\loc(\bbr_+\times\bbr^d)\cap L^2(\bbr_+\times\bbr^d)$ but $g-\rho*g\notin L^1(\bbr_+\times\bbr^d)$, and $g-\rho*g$ is non-negative or non-positive for all $(t,x)\in \bbr_+\times\bbr^d$, then the process $X$ in \eqref{twosidedX} is well defined, has a finite second moment and is long-range dependent in the sense of \eqref{longmem}.
\eenu
\eprop

\begin{Example}[Long-range dependence by temporal regular variation]\label{lrd-ex1} ~
	
	\noindent
	We consider a measure $\mu$ of the form $\mu(\dd t,\dd x) = k(t)\,\dd t\,\delta_{0,\bbr^d}(\dd x)$ with some  $k\in L^1_\loc(\bbr_+)$. By Lemma~\ref{resLemma}, the resolvent of $\mu$ has the form $\rho(\dd t,\dd x) = r(t)\,\dd t\,\delta_{0,\bbr^d}(\dd x)$ for some $r\in L^1_\loc(\bbr_+)$, which is the temporal resolvent of the function $k$ (i.e., we have $r+k=r*k$ where $*$ stands for convolution on $\bbr_+$). Now suppose that the function $k$ satisfies the following assumptions:
	\bit
	\item $k(t)=-t^{-\al}L(t)$ for all $t\in\bbr_+$, some $\al\in(0,1/2)$ and some function $L\colon\bbr_+\to(0,\infty)$ that is slowly varying at infinity.
	\item $k$ is differentiable with a continuous strictly positive derivative that belongs to $L^1(\bbr_+)$.
	\item The function $t\mapsto \log(-k(t))$ is convex in $t$.
	\eit
	Then, by Theorem~3.2 of \cite{Appleby11},
	\beq\label{regvar} \lim_{t\to\infty} \left(1-\int_0^t r(s)\,\dd s \right)t^{1-\al}L(t)=\frac{\sin(\al\pi)}{\pi}. \eeq 
	If now $g(t,x)=g_0(x)$ for some non-negative (or non-positive) $g_0\in L^1(\bbr^d)\cap L^2(\bbr^d)$, then 
	\[ (g-\rho*g)(t,x)=g_0(x)\left(1-\int_0^t r(s)\,\dd s\right),\quad (t,x)\in\bbr_+\times\bbr^d,\]
	is, because of \eqref{regvar} and Corollary~8.8 of Chapter~9 in \cite{Gripenberg90}, non-negative (or non-positive), belongs to $L^1_\loc(\bbr_+\times\bbr^d)\cap L^2(\bbr_+\times\bbr^d)$, but not to $L^1(\bbr_+\times\bbr^d)$. Hence we conclude from Proposition~\ref{SRD} that the resulting stationary process in \eqref{twosidedX} has long-range dependence. One possible choice of $k$ is $k(t)=-1/(\al(1+t)^\al)$ with $\al\in(0,1/2)$, cf. Example~3.5 in \cite{Appleby11}. 	
	\halmos
\end{Example}
In the previous example, the non-integrability of the resolvent measure is essentially due to the regular variation of the function $k$. In the next example, long-range dependence arises through a drift measure of the form $\mu=-\la \Leb_{\bbr_+}\otimes m$ with some $\la>0$ and $m\in M(\bbr^d)$. 
\begin{Example}[Long-range dependence by spatial regular variation]\label{lrd-ex2} ~
	
	\noindent
	Consider the measure $\mu(\dd t,\dd x) = -\dd t\,f(x)\,\dd x$ with $f(x)=1/(\pi(1+x^2))$. In Example~\ref{Exres}(2) the resolvent measure is found to have the Lebesgue density $r(t,x)=(2\pi x)^{-1} G(t,x)$ where $G$ is the function given in \eqref{funcG}.
	We use the software package Mathematica to check that $r\notin L^1(\bbr_+\times\bbr)$.
	\begin{comment}
	To see that $r\notin L^1(\bbr_+\times\bbr)$, we use the software package Mathematica to show that for all $x\neq0$ the definite integral
	\[ \int_0^T r(t,x)\,\dd t = \frac{\ii}{2\pi x}\Big(T^{-\ii x}(\Ga(1+\ii x)-\Ga(1+\ii x,T))-T^{\ii x} (\Ga(1-\ii x)-\Ga(1-\ii x, T))\Big)  \]
	does not converge as $T\to \infty$. While
	\[ \frac{\ii}{2\pi x}(T^{\ii x}\Ga(1-\ii x,T)-T^{-\ii x}\Ga(1+\ii x,T))\to 0,\quad T\to \infty,\]
	the map $T\mapsto \ii/(2\pi x)(T^{-\ii x}\Ga(1+\ii x)-T^{\ii x}\Ga(1-\ii x))$ is oscillatory: it takes the value 
	\[ (-1)^n(\pi x \sinh(\pi x))^{-1/2}\quad \text{at} \quad T=\ee^{-\pi (2n-1)/(2x)}(\Ga(1-\ii x)/\Ga(1+\ii x))^{\ii/(2x)},\quad n\in\bbn\]
	(these are real-valued despite their appearance).
	This proves that $r\notin L^1(\bbr_+\times\bbr)$. 
	\end{comment}
	However, we do have that $r\in L^2(\bbr_+\times\bbr)$. Indeed, if 
	\[ \tilde f(x):=\frac{1}{\sqrt{2\pi}} \int_\bbr \ee^{-\ii x z} f(z)\,\dd z,\quad x\in\bbr,  \]
	denotes the Fourier transform of a function $f\in L^2(\bbr)$, we can use the fact that $\widetilde{f^{*n}} = \tilde f^n$ and $\tilde f(x)=(2\pi)^{-1/2}\ee^{-|x|}$ for $f(x)=1/(\pi(1+x^2))$ to calculate the Fourier transform of $r(t,\cdot)$ for fixed $t\in\bbr_+$:
	\[ \tilde r(t,x)= \sum_{n=1}^\infty \frac{(-t)^{n-1}}{(n-1)!} (\tilde f(x))^n = \frac{1}{\sqrt{2\pi}}\ee^{-|x|}\ee^{-t(2\pi)^{-1/2}\ee^{-|x|}},\quad (t,x)\in\bbr_+\times\bbr. \]
	Since Fourier transformation is unitary on $L^2(\bbr)$, by Plancherel's theorem, we have
	\[ \int_{\bbr_+} \int_\bbr r(t,x)^2\,\dd (t,x) = \int_{\bbr_+}\int_{\bbr} \tilde r(t,x)^2\,\dd (t,x)=\frac{1}{2\pi} \int_{\bbr_+} \frac{\pi-\ee^{-\sqrt{2/\pi}t}(\pi + \sqrt{2\pi} t)}{t^2} \,\dd t = \frac{1}{\sqrt{2\pi}}. \]
	\halmos
\end{Example}

\section{Path properties}\label{pathprop}
The classical OU process \eqref{OUsol} has nice path properties: if $W$ is Gaussian, the process has a continuous version; if $W$ has jumps, the solution has a càdlàg version. In contrast, the notion of solution in Theorem~\ref{localsol} basically says nothing about the paths of a VOU process (apart from being measurable and elements of the set $\mathcal{L}(\mu)$ in \eqref{Lmu}). The goal of this section is to fill in this gap, at least partially, and to prove the existence of versions with nice regularity for the VOU process. In the presence of jumps, it turns out that temporo-spatial path properties are much harder to establish than for processes indexed by time. But before discussing this in detail, we first consider the case where the driving noise is Gaussian and mild conditions already ensure the existence of a Hölder continuous version.
%In this section we analyze path properties of our solution process. Our goal is to prove the existence of a version with càdlàg sample paths. However, there is no canonical notion of càdlàg in space and time. It turns out that, depending on the regularity of the noise propagation function $g$ in the VOU equation \eqref{VOU}, different kinds of temporo-spatial càdlàg properties emerge. Nonetheless, we start with a useful criterion for continuous sample paths in the case where the Lévy basis is a Gaussian basis with drift. 
%In the following let $W_{\loc}^{k,p}(\bbr_+\times\bbr^d)$ denote the collection of all functions $\bbr_+\times\bbr^d\to\bbr$ whose restrictions to $[0,T]\times\bbr^d$ belong to the Sobolev space $W^{k,p}([0,T]\times\bbr^d)$ for all $T\in\bbr_+$ (see Section~5.2 of \cite{Evans98} for a definition of the Sobolev spaces). 
%i.e. the mixed partial derivatives up to order $k$ exist in the weak sense and also lie in $L^p([0,T]\times\bbr^d)$.
\begin{Theorem}\label{contversion}
	Suppose that the conditions of Theorem~\ref{localsol} are satisfied with $V\equiv0$, $\nu\equiv0$ and that $X$ is the solution to \eqref{VOU} given in \eqref{solformula}. Further assume:
	\begin{itemize}
		\item There exists an exponent $u>0$ such that for every $T>0$ there is a non-negative constant $C_T$ and we have
		\beq\label{translation} \int_0^T\int_{\bbr^d}\left| g(s,y)- g(s+\tau,y+\xi)\right|^2\,\dd(s,y)\leq C_{T}|(\tau,\xi)|^u \eeq
		whenever $|(\tau,\xi)|$ is sufficiently small.
		\item For some $p>1$  we have for every $T>0$ that
		\beq\label{intp} \int_0^T\left(\int_{\bbr^d} g(s,y)^2\,\dd y\right)^p\,\dd s<\infty.\eeq 
		%\item there are $v>0$ and $K>0$ such that $\left(\int_t^{t+\tau}\int_{\bbr^d} g(s,y)\,\dd (s,y)\right)^2\leq K\tau^v$ for every $t>0$	and sufficiently small $\tau>0$.	
	\end{itemize} 
	%$\tilde g:=g-\rho*g\in L_{\loc}^p(\bbr_+\times\bbr^d)\cap W_{\loc}^{1,q}(\bbr_+\times\bbr^d)$ for some $p,q \in (1,\infty)$ with $p^{-1}+q^{-1}=1$. 
	Then the process $X$ has a version which is locally Hölder continuous with any exponent in $(0,\frac{p-1}{4p}\wedge\frac{u}{4})$.
\end{Theorem}

As usual for Gaussian processes, continuity of sample paths can be established under weaker conditions than those formulated in Theorem~\ref{contversion}. Given that $g\in L^2_\loc(\bbr_+\times\bbr^d)$ already implies that the left-hand side of \eqref{translation} converges to $0$ as $|(\tau,\xi)|\to0$ and \eqref{intp} holds with $p=1$, the assumptions of Theorem~\ref{contversion} are reasonably general for practical purposes, so we do not pursue this direction further and only refer to \cite{Adler81}

If the noise does feature jumps, we are not able to construct path properties in general. Informally speaking, if the noise propagation function $g$ is too irregular, it is unclear how a jump at a certain time and location affects the process at other times and locations. However, the situation is different if $g$ is smooth enough. In this case, we are able to establish versions with the following regularity property.
\begin{Definition}\label{deftempcadlag}
	A function $\Psi\colon\bbr_+\times\bbr^d\to\bbr$ is  \emph{t-càdlàg} if for every $(t,x)\in\bbr_+\times\bbr^d$,
	\[ \lim_{\substack{(\tilde t,\tilde x)\to(t, x)\\ \tilde t\geq t}}\Psi(\tilde t,\tilde x)=\Psi(t,x)\quad\text{and}\quad \lim_{\substack{(\tilde t,\tilde x)\to(t, x)\\ \tilde t < t}}\Psi(\tilde t,\tilde x)\quad\text{exists.} \]
	\halmos
\end{Definition}
In the following $\pd^\ga g$ denotes the partial derivative $\pd_t^{\ga_0}\pd_{x_1}^{\ga_1}\cdots\pd_{x_d}^{\ga_{d}}g$ for a function $g\colon\bbr_+\times\bbr^d\to\bbr$ and a multi-index $\ga=(\ga_0,\ldots,\ga_{d})\in\bbn_0^{d+1}$.
\begin{Theorem}\label{tempcadlag}
	Suppose that the conditions of Theorem~\ref{localsol} are satisfied with $V\equiv0$ and that $X$ is the solution to \eqref{VOU} given in \eqref{solformula}. We assume that $g$ is $(d+1)$-times continuously differentiable on $\bbr_+\times\bbr^d$ such that for all multi-indices $\ga\in\{0,1\}^{d+1}$ the partial derivative $\pd^\ga g$ belongs to $L^1_\loc(\bbr_+\times\bbr^d)$. 
	
	If $\al<1$, we further assume that there is a non-negative decreasing function $G\colon \bbr_+\to\bbr_+$ such that $G(|x|)$ belongs to $L^\al(\bbr^d)$ and we have $| g(t,x)|\leq C_T G(|x|)$ for all $T\geq 0$ and $(t,x)\in[0,T]\times\bbr^d$, where $C_T$ is a non-negative constant depending on $T$.
	
	Then the process $X$ has a t-càdlàg version. This version is continuous if $g$ additionally satisfies $g(0,x)=0$ for all $x\in\bbr^d$.
	%\begin{align*}
	%	&\int_0^t\int_\bbr \left[\int_0^x |(\pd_2g)(0,z-y)|\,\dd z\right]^\al\,\dd (s,y)<\infty,\\
	%	&\int_0^t\int_\bbr \left[\int_s^t |(\pd_1g)(r-s,-y)|\,\dd r\right]^\al\,\dd (s,y)<\infty,\\
	%	&\int_0^t\int_\bbr \left[\int_s^t\int_0^x |(\pd_1\pd_2g)(r-s,z-y)|\,\dd (r,z)\right]^\al\,\dd (s,y)<\infty
	%\end{align*}
	%If $g$ additionally satisfies $g(0,x)=0$ for all $x\in\bbr^d$, this version is continuous.
\end{Theorem}

The solution may display a fundamentally different path behavior if the underlying noise propagation function is not smooth. Here we only have results for the spatial dimension $1$ and $g(t,x)=\bone_{\{|x|\leq ct\}}h(t,x)$ with some $c>0$ and some smooth function $h$. This choice for $g$ is motivated by Example~\ref{Ex2} and, as we shall see in the proofs, enables us to utilize maximal inequalities of multi-parameter martingales as in \cite{Walsh79}. 

In order to state our result, we introduce a partial order $\lpo$ on $\bbr_+\times\bbr$ by setting $(t, x)\lpo(\tilde t,\tilde x)$ if $t\leq\tilde t$ and $|\tilde x-x|\leq c(\tilde t-t)$. As usual, we write $(t, x)\prec(\tilde t,\tilde x)$ if $(t, x)\lpo(\tilde t,\tilde x)$ and $(t, x)\neq(\tilde t,\tilde x)$. The following temporo-spatial càdlàg property is weaker than the t-càdlàg property of Definition~\ref{deftempcadlag}.
\begin{Definition}\label{defpocadlag}
	A function $\Psi\colon\bbr_+\times\bbr^d\to\bbr$ is  \emph{$\lpo$-càdlàg} if for every $(t,x)\in\bbr_+\times\bbr^d$,
	\[ \lim_{\substack{(\tilde t,\tilde x)\to(t, x)\\(t, x)\lpo(\tilde t,\tilde x)}}\Psi(\tilde t,\tilde x)=\Psi(t,x)\quad\text{and}\quad\lim_{\substack{(\tilde t,\tilde x)\to(t, x)\\(t, x)\succ(\tilde t,\tilde x)}}\Psi(\tilde t,\tilde x)\quad\text{exists.} \]
	\halmos
\end{Definition}
\begin{Theorem}\label{pocadlag}
	Suppose that the conditions of Theorem~\ref{localsol} are satisfied with $d=1$, $V\equiv0$ and that $X$ is the solution to \eqref{VOU} given in \eqref{solformula}. Further let $h\in L^\al_\loc(\bbr_+\times\bbr)\cap L^\beta_\loc(\bbr_+\times\bbr)$ be twice continuously differentiable with partial derivatives $\partial_t h(t,x), \partial_x h(t,x)$ and $\partial_t \partial_x h(t,x)$ in $L^1_\loc(\bbr_+\times\bbr)$. If the noise propagation function takes the form $g(t,x)=\bone_{\{|x|\leq ct\}}h(t,x)$ for some $c>0$, then the process $X$ has a $\lpo$-càdlàg version.
\end{Theorem}

\begin{Remark}\label{remcadlag}
	\begin{enumerate}
		\item The conditions of Theorem~\ref{tempcadlag} and Theorem~\ref{pocadlag} imply in particular that $g$ is bounded. It is important to notice that the assertions of these theorems are false when $g$ has a singularity at, for example, the origin.
		\item Under the assumptions of Theorem~\ref{pocadlag} and with virtually no change of its proof, we even have a version of the process $X$ which is not only $\lpo$-càdlàg but also has limits from the flanks, that is, both limits
		\[ \lim_{\substack{(\tilde t,\tilde x)\to(t, x)\\\tilde x>x,\, c|t-\tilde t|\leq\tilde x-x}}X(\tilde t,\tilde x)\quad\text{and}\quad\lim_{\substack{(\tilde t,\tilde x)\to(t, x)\\\tilde x<x,\, c|t-\tilde t|\leq x-\tilde x}}X(\tilde t,\tilde x)\quad\text{exist.} \] 
		\item There exist other notions of càdlàg sample paths for multi-parameter stochastic processes, see e.g. \cite{Neuhaus71}. In contrast to the definition in that reference, our Definitions~\ref{deftempcadlag} and \ref{defpocadlag} take into account that time has a natural direction, while space has none.
		\item It suffices for Theorem~\ref{tempcadlag} (resp. Theorem~\ref{pocadlag}) in dimension $d=1$ that, instead of being twice continuously differentiable, $g$ (resp. $h$) is continuous on $\bbr_+\times\bbr$ and that there exist $\beta\in[1,2]$ (resp. $\beta\in(1,2]$) and functions $g_1, g_2, g_{12}\in L^\beta_\loc(\bbr_+\times\bbr)$ (resp. $h_1, h_2, h_{12}\in L^\beta_\loc(\bbr_+\times\bbr)$) such that $\int_\bbr |z|^\beta\bone_{\{|z|\leq1\}}\,\nu(\dd z)<\infty$, Equation~\eqref{help5} (resp. Equation~\eqref{help2}) holds and, if $\al=1$ in Theorem~\ref{tempcadlag}, that also $g_1, g_2, g_{12}\in L^1_\loc(\bbr_+\times\bbr)$. For higher dimensions in Theorem~\ref{tempcadlag}, this comment applies analogously.
	\end{enumerate}		
	\halmos
\end{Remark}

We apply the derived theorems to the VOU model considered in Examples~\ref{Ex1} and \ref{Ex2}.

\begin{Example}[Path properties for Examples~\ref{Ex1} and \ref{Ex2}]\label{expath}~
	\begin{enumerate}
				
	\item By induction we can show that the partial derivatives of the noise propagation function $g(t,x)=\ee^{-\la'|x|}$ can be written as 
	\[ \partial_{x_1}\ldots\partial_{x_n} g(t,x)=\sum_{j=n}^{2n-1}c_j\ee^{-\la'|x|}\frac{x_1\ldots x_n}{|x|^j}, \]
	for all $n\leq d$ and constants $c_j$ independent of $(t,x)$. As a consequence, every partial derivative $\pd^\ga g$, where $\ga$ is some multi-index in $\{0,1\}^{d+1}$, belongs to $L^p_\loc(\bbr_+\times\bbr^d)$ for every $0<p<\frac{d}{d-1}$. For $\al<1$ we may choose the function $G$ from Theorem~\ref{tempcadlag} as $G(u)=\ee^{-\la'u}$. Therefore, as soon as the L\'evy measure $\nu$ of $\La$ satisfies \eqref{momLa}
	for some $\al>0$ and $\beta\in(0,\frac{d}{d-1})$, Theorem~\ref{tempcadlag} and Remark~\ref{remcadlag} imply that the unique solution $X$ in \eqref{eqsolmodel1} has a t-càdlàg version.

	\item Since the partial derivatives $\partial_t h(t,x), \partial_x h(t,x)$ and $\partial_t \partial_x h(t,x)$ of the function $h(t,x)=\ee^{-\la'|x|}$ satisfy $\partial_t h(t,x)=0, \partial_x h(t,x)=-\la'\ee^{-\la'|x|}\frac{x}{|x|}$ and $\partial_t \partial_x h(t,x)=0$ and lie in $L^p_\loc(\bbr_+\times\bbr^d)$ for every $p>0$, Theorem~\ref{pocadlag} and Remark~\ref{remcadlag} apply to Example~\ref{Ex2}. Hence, we obtain a $\lpo$-càdlàg version for the process $X$ in \eqref{eqsolmodel2}.
	\halmos
	
	\end{enumerate}	
\end{Example}

\section{Proofs}\label{sectproofs}
For the rest of this paper, $C_T$ denotes a real constant which may depend on $T\geq0$ and change its value from line to line. 
\subsection{Proofs for Section~\ref{soltheory}}
For the proofs of Theorems~\ref{localsol} and \ref{stationarysol} we must guarantee that the stochastic convolution of $g$ with $\La$ is well defined on $[0,T]\times\bbr^d$ or on the whole $\bbr_+\times\bbr^d$, respectively. The conditions listed in Proposition~\ref{RR} are necessary and sufficient, but may be too complicated to verify in general. The following lemma provides some simpler sufficient criteria. 
\blem\label{integrability-simple}
Suppose that $\La$ is a homogeneous L\'evy basis on $I\times\bbr^d$ with characteristics $(b,\si^2,\nu)$ and that $g\colon I\times\bbr^d\to\bbr$ is a measurable function.
\benu

\item Condition (1) of Proposition~\ref{RR} is satisfied if there exist measurable functions $k, K\colon \bbr_+ \to \bbr_+$ such that
\begin{align}
	\label{infassym} \left|b+\int_\bbr z\bone_{\{|z|\in(1,A]\}} \,\nu(\dd z)\right| &= \Omicron(K(A)),\quad A\to\infty, \\
	\label{0assym} \left|b-\int_\bbr  z\bone_{\{|z|\in (a,1]\}}\,\nu(\dd z)\right| &= \Omicron(k(a)),\quad a\to0, \quad\text{and}
\end{align}
\beq\label{kK} \int_I \int_{\bbr^d}  |g(t,x)| \Big( k(|g(t,x)|^{-1}) \bone_{\{|g(t,x)|>1\}} + K(|g(t,x)|^{-1})\bone_{\{|g(t,x)|\leq1\}} \Big)\,\dd(t,x) < \infty.
\eeq

\item Condition (3) of Proposition~\ref{RR} is satisfied if there is an increasing function $h\colon\bbr_+\to \bbr_+$ such that for some constant $C\in\bbr_+$ we have $u^2 h(|x|)\leq C h(u|x|)$ for all $u\in[0,1]$ and $x\in\bbr$, and
\begin{equation*}
	%\label{hassump} 
	\int_I \int_{\bbr^d} h(|g(t,x)|) \,\dd(t,x) < \infty \quad\text{and}\quad \int_\bbr h(|z|^{-1})^{-1}\,\nu(\dd z) < \infty. 	
\end{equation*}

\eenu
\elem

\bpr%[of Lemma~\ref{integrability-simple}]
(1) The left-hand side of condition (1) of Proposition~\ref{RR} is bounded by 
\begin{align*} &~\int_I\int_{\bbr^d}  |g(t,x)|\left|b-\int_\bbr z\bone_{\{|z|\in (|g(t,x)|^{-1},1]\}} \,\nu(\dd z)\right|\bone_{\{|g(t,x)|>1\}}\,\dd (t,x)\\
	&~\qquad +\int_I\int_{\bbr^d}  |g(t,x)| \left|b+\int_\bbr z \bone_{\{|z|\in(1,|g(t,x)|^{-1}]\}} \,\nu(\dd z)\right| \bone_{\{|g(t,x)|\leq1\}}\,\dd (t,x),%\\
	%\leq&~ \iint_{|g(t,x)|>1} |g(t,x)|k(|g(t,x)|^{-1})\,\dd (t,x) + \iint_{|g(t,x)|\leq1} |g(t,x)|K(|g(t,x)|^{-1})\,\dd (t,x) < \infty.
\end{align*}
which, because of \eqref{infassym} and \eqref{0assym}, is in turn bounded by \eqref{kK}.

\vspace{0.5\baselineskip}\noindent (2) We divide the integral term in condition (3) of Proposition~\ref{RR} into
\[ J_1+J_2:=(\Leb_{I\times\bbr^d}\otimes\nu)(\{|zg(t,x)|>1\}) + \int_{I}\int_{\bbr^d}\int_\bbr  |z g(t,x)|^2\bone_{\{|zg(t,x)|\leq 1\}}\,\nu(\dd z)\,\dd(t,x). \]
For $J_1$ we can now use Markov's inequality to obtain 
\[ J_1\leq \int_\bbr h(|z|^{-1})^{-1}\,\nu(\dd z) \int_I\int_{\bbr^d} h(|g(t,x)|)\,\dd (t,x) <\infty,\]
while for $J_2$ the assumption that $|zg(t,x)|^2\leq C h(|g(t,x)|)h^{-1}(|z|^{-1})$ on $\{|zg(t,x)|\leq1\}$ implies
\[ J_2 \leq C \int_I \int_{\bbr^d} h(|g(t,x)|)\,\dd (t,x) \int_\bbr h(|z|^{-1})^{-1}\,\nu(\dd z) < \infty.\] \halmos
\epr

Possible choices for the functions $h,k$ and $K$ are $h(x):=x^q\bone_{[0,1]}(x) +x^p\bone_{(1,\infty)}(x)$, $k(a):=a^{1-p}$ and $K(A):=A^{1-q}$ with some $0<p,q\leq2$.

\bpr[of Theorem~\ref{localsol}] %Without loss of generality, we may assume that $\al\in(0,1]$ and $\beta\in[1,2]$ because all conditions in Theorem~\ref{localsol} continue to hold if $\al$ is replaced by $\al\wedge 1$ and $\beta$ by $\beta\vee1$.

We wish to apply part (3) of Proposition~\ref{thmdetsolution} to the stochastic forcing function defined by
\[ F(t,x):=V(t,x)+\int_0^t\int_{\bbr^d} g(t-s,x-y)\,\La(\dd s,\dd y),\quad(t,x)\in\bbr_+\times\bbr^d. \]
This would yield the existence and uniqueness statement of Theorem~\ref{localsol}. Since $V$ has paths in $\mathcal{F}(\mu)$ by hypotheses, it suffices to prove that 
\[ Y(t,x):=\int_0^t\int_{\bbr^d} g(t-s,x-y)\,\La(\dd s,\dd y)\] is well defined for all $(t,x)\in\bbr_+\times\bbr^d$ and that $Y$ has a version with paths in $\mathcal{F}(\mu)$. The existence of the stochastic convolution is equivalent to the integrability of $g$ with respect to $\La$. But this follows from Lemma~\ref{integrability-simple} with $h(x):=x^\al\bone_{[0,1]}(x) +x^\beta\bone_{(1,\infty)}(x)$, $k(a):=a^{1-\beta}$ and $K(A):=A^{1-\al}$ because we have $g\in L^\al_\loc(\bbr_+\times\bbr^d)\cap L^\beta_\loc(\bbr_+\times\bbr^d)$ on the one hand, and \eqref{momLa} on the other hand. Note that the latter also implies \eqref{infassym} and \eqref{0assym} with our choices of $k$ and $K$. Furthermore, from the well-definedness of $Y$ we can already deduce the existence of a measurable version, see Theorem~1 of \cite{Lebedev96}.

Next, we prove that this measurable version of $Y$ belongs to $\mathcal{F}(\mu)$ almost surely. We begin with the case $\al=1$ and notice that applying the Jensen and the Burkholder-Davis-Gundy inequalities yields
\begin{align} \bbe[|Y(t,x)|]\leq&~ |b|\int_0^t \int_{\bbr^d} |g(s,y)|\,\dd (s,y) + \left( \int_0^t\int_{\bbr^d} \si^2|g(s,y)|^2\,\dd (s,y) \right)^{1/2} \nonumber \\
&~+ C_T\bbe\left[ \left( \int_0^t\int_{\bbr^d} \int_{\bbr} |g(t-s,x-y)z|^2 \bone_{\{|z|\leq 1 \}}\,\pf(\dd s,\dd y,\dd z) \right)^{\beta/2} \right]^{1/\beta} \nonumber\\
&~+ \bbe\left[\int_0^t\int_{\bbr^d}\int_\bbr |g(t-s,x-y)z|\bone_{\{|z|> 1 \}} \,\pf(\dd s,\dd y,\dd z)  \right]\nonumber\\
\leq &~ |b|\int_0^t \int_{\bbr^d} |g(s,y)|\,\dd (s,y) + \left( \int_0^t\int_{\bbr^d} \si^2|g(s,y)|^2\,\dd (s,y) \right)^{1/2}\nonumber \nonumber\\
&~+ C_T\left( \int_0^t\int_{\bbr^d}  |g(s,y)|^\beta \,\dd(s,y) \int_{\bbr} |z|^\beta \bone_{\{|z|\leq 1 \}}\,\nu(\dd z) \right)^{1/\beta}\nonumber \\
&~+ \int_0^t\int_{\bbr^d} |g(s,y)|\,\dd (s,y) \int_{\bbr} |z|\bone_{\{|z|> 1 \}} \,\nu(\dd z).\label{help0} \end{align}
Therefore, the function $(t,x)\mapsto \bbe[|Y(t,x)|]$ belongs to $L^\infty_\loc(\bbr_+\times\bbr^d)$, so Lemma~\ref{prop2} and Fubini's theorem imply that $Y$ has paths in $\mathcal{F}(\mu)$ almost surely.

If $\al\in(0,1)$, it is enough to prove $Y^1\in \mathcal{F}(\mu)$ almost surely where $Y^1$ is defined in the same way as $Y$ but with $\La$ replaced by its large jumps part $\La^1(\dd t,\dd x) := \int_\bbr z\bone_{\{|z|> 1 \}} \,\pf(\dd t,\dd x,\dd z)$. For the convolution of $g$ with $\La-\La^1$ the arguments as in the case $\al=1$ would apply. Letting $\eta\in\{|\mu|, |\rho|, |\mu|*|\rho|, |\mu|*|\rho|^{*2}\}$, we prove that $\eta*|Y^1|$ exists almost everywhere. Since the realizations of $\La^1$ are measures on $\bbr_+\times\bbr^d$, Fubini's theorem yields
\begin{align*} (\eta*|Y^1|)(t,x) =&~ \int_0^t\int_{\bbr^d} |Y^1(t-s,x-y)|\,\eta(\dd s,\dd y)\\ 
	\leq&~ \int_0^t\int_{\bbr^d} \int_0^{t-s}\int_{\bbr^d} \int_\bbr |g(t-s-r,x-y-w)z|\bone_{\{|z|> 1 \}} \,\pf(\dd r,\dd w,\dd z)\,\eta(\dd s,\dd y)\\
	=&~ \int_0^t \int_{\bbr^d}\int_{\bbr} (\eta*|g|)(t-s,x-y)|z|\bone_{\{|z|> 1 \}}\,\pf(\dd s,\dd y,\dd z). 
\end{align*}
Now raising the last inequality to the power $\al$, moving the exponent into the integral and taking expectation result in
\begin{align*} \bbe\big[\big((\eta*|Y^1|)(t,x)\big)^\al\big] \leq&~ \int_0^t\int_{\bbr^d} (\eta*|g|)^\al(s,y)\,\dd (s,y) \int_\bbr |z|^\al\bone_{\{|z|> 1 \}}\,\nu(\dd z)\\
=&~\int_0^t\int_{\bbr^d} \big(\vp(y)(\eta*|g|)(s,y)\big)^\al \vp^{-\al}(y)\,\dd (s,y) \int_\bbr |z|^\al\bone_{\{|z|> 1 \}}\,\nu(\dd z). \end{align*}
Since $\vp^{-\al} \in L^1(\bbr^d)$, the assertion that $Y^1\in \mathcal{F}(\mu)$ almost surely is proved once we can show that $\vp(\eta*|g|) \in L^\infty_\loc(\bbr_+\times\bbr^d)$. For $\eta=|\mu|$ this holds by assumption. For $\eta=|\rho|$ we use the fact that $\rho=\rho*\mu-\mu$, so we only need to prove $\vp(|\rho|*(|\mu|*|g|)) \in L^\infty_\loc(\bbr_+\times\bbr^d)$. Notice that $\vp\mu \in M_\loc(\bbr_+\times\bbr^d)$ implies $\vp\rho\in M_\loc(\bbr_+\times\bbr^d)$ by the same arguments as in Proposition~\ref{thmdetsolution}. More precisely, one has to work in the weighted measure spaces $M([0,T]\times\bbr^d;\vp)$, $T\in\bbr_+$, which consist of all signed complete Borel measures $\mu$ such that $\vp\mu\in M([0,T]\times\bbr^d)$ and are complete with the weighted total variation norm $\|\cdot\|_\vp := \|\vp \cdot\|$ (the proof of the temporal analogue, Theorem~4.3.4 of \cite{Gripenberg90}, can be extended to the temporo-spatial setting in a straightforward manner). Thus, the hypothesis that $\vp(|\mu|*|g|) \in L^\infty_\loc(\bbr_+\times\bbr^d)$ yields $\vp(|\rho|*(|\mu|*|g|)) \in L^\infty_\loc(\bbr_+\times\bbr^d)$ as well (like before, one can extend Theorem~4.3.5 in \cite{Gripenberg90}). Finally, for $\eta=|\mu|*|\rho|^{*2}$ the same arguments apply because we have already established $\vp(|\rho|*(|\mu|*|g|)) \in L^\infty_\loc(\bbr_+\times\bbr^d)$.

It remains to demonstrate that \eqref{solformula} defines a version of the solution to \eqref{VOU}. To this end, we first observe that the solution in Proposition~\ref{thmdetsolution} takes the form
\begin{align*} X(t,x)=&~ V(t,x) - (\rho*V)(t,x) + \int_0^t \int_{\bbr^d} g(t-s,x-y)\,\La(\dd s,\dd y)\\
&~ - \int_0^t\int_{\bbr^d} \int_0^{t-s}\int_{\bbr^d} g(t-s-r,x-y-w)\,\La(\dd r,\dd w)\,\rho(\dd s,\dd y). \end{align*}
Formula \eqref{solformula} immediately follows if we can interchange the integrals with respect to $\La$ and $\rho$, that is, if we can apply a stochastic Fubini theorem. For the large jumps part $\La^1$ of $\La$ the ordinary Fubini theorem is sufficient because the realizations of $\La^1$ are true measures and integrability has already been shown in the proof for $|\rho|*|Y^1|$ above. For the remaining part $Y^2:=Y-Y^1$ Theorem~2 in \cite{Lebedev96} is applicable because, by the same reasoning as in \eqref{help0}, we have
\[ \int_0^t\int_{\bbr^d} \bbe\left[ \left(\int_0^{t-s}\int_{\bbr^d}\int_\bbr  |g(t-s-r,x-y-w)z|^2\bone_{\{|z| \leq1 \}}\,\pf(\dd r,\dd w,\dd z)\right)^{1/2} \right]|\rho|(\dd s,\dd y)<\infty. \]
\halmos
\epr

\bpr[of Theorem~\ref{stationarysol}] Our first observation is that $X_\infty$ is well defined because $g-\rho*g$ is integrable with respect to $\La$ on $\bbr_+\times\bbr^d$. This in turn is a consequence of assumptions (2) and (3) together with Lemma~\ref{integrability-simple} (and of course, that $\beta=2$ if $\si^2\neq0$). Next, regarding the convergence statement, it suffices by Slutsky's theorem and hypothesis (1) to prove the convergence of the finite-dimensional distributions of the stochastic convolution in \eqref{solformula} when time tends to infinity. For one spatial point $x\in\bbr^d$, the claim readily follows from
\begin{align*} \int_0^t\int_{\bbr^d} (g-\rho*g)(t-s,x-y)\,\La(\dd s, \dd y) \stackrel{\dd}{=}&~ \int_0^t\int_{\bbr^d} (g-\rho*g)(s,x-y)\,\La(t-\dd s,\dd y)\\ \stackrel{\dd}{=}&~ \int_0^t\int_{\bbr^d} (g-\rho*g)(s,x-y)\,\La(\dd s,\dd y)\\ \stackrel{\dd}{\longrightarrow} &~\int_0^\infty\int_{\bbr^d} (g-\rho*g)(s,x-y)\,\La(\dd s,\dd y). \end{align*}
The $n$-dimensional case can be treated completely analogously. For the second part of the theorem, we suppose that also $g$ is integrable with respect to $\La$ and that $V$ is a process independent of $\La$ and with the same finite-dimensional distributions as the process given in \eqref{Vdist}. Then by the stochastic Fubini theorem (cf. the proof of Theorem~\ref{localsol}) we obtain
\begin{align*} &~V(t,x)-(\rho*V)(t,x)\\
&~\qquad \eqd \int_{-\infty}^0 \int_{\bbr^d} g(t-s,x-y)\,\La(\dd s,\dd y) \\
&~\qquad\quad- \int_0^t\int_{\bbr^d} \int_{-\infty}^0 \int_{\bbr^d} g(t-s-r,x-y-w)\,\La(\dd r,\dd w)\,\rho(\dd s,\dd y)\\
&~\qquad \eqd\int_{-\infty}^0 \int_{\bbr^d} g(t-s,x-y)\,\La(\dd s,\dd y) - \int_{-\infty}^0 \int_{\bbr^d} (\rho*g)(t-s,x-y)\,\La(\dd s,\dd y).  \end{align*}
Again, the reader can convince herself that the previous calculations also apply to $n$ time and space points. The strict stationarity of $X$ is now a consequence of that of the process \eqref{twosidedX}.
\halmos
\epr

\bpr[of Lemma~\ref{statcondLa}] 
(1) We have already used this tacitly in the proof of Theorem~\ref{localsol}. In fact, if $\al\in (0,1]$, then
\[ \left| b+\int_\bbr z \bone_{\{ |z|\in(1,A] \}} \,\nu(\dd z)\right|\leq |b|+\int_\bbr |z|^\al |z|^{1-\al}\bone_{\{ |z|\in(1,A] \}} \,\nu(\dd z) \leq |b|+A^{1-\al}\int_\bbr |z|^\al\bone_{\{ |z|>1 \}}\,\nu(\dd z), \]
while for $\beta\in[1,2]$ we have 
\[ \left| b-\int_\bbr z \bone_{\{ |z|\in (a,1] \}} \,\nu(\dd z)\right|\leq |b|+\int_\bbr |z|^\beta |z|^{1-\beta}\bone_{\{ |z|\in(a,1] \}} \,\nu(\dd z) \leq |b|+a^{1-\beta}\int_\bbr |z|^\beta\bone_{\{ |z|\leq 1 \}}\,\nu(\dd z). \]

\vspace{0.5\baselineskip}\noindent (2) If $\al\in(1,2]$ and $b_1=0$, then we have for $A\geq1$ that
\begin{align*} \left| b+\int_\bbr z \bone_{\{ |z|\in(1,A] \}} \,\nu(\dd z)\right|\leq&~ \int_\bbr |z|\bone_{\{ |z|\in(A,\infty) \}} \,\nu(\dd z) =  \int_\bbr |z|^\al |z|^{1-\al}\bone_{\{ |z|\in(A,\infty) \}} \,\nu(\dd z)\\
\leq&~ A^{1-\al}\int_\bbr |z|^\al\bone_{\{ |z|>1 \}}\,\nu(\dd z).\end{align*}
If $\beta\in(0,1)$ and $b_0=0$, then 
\begin{align*}
\left| b-\int_\bbr z \bone_{\{ |z|\in(a,1] \}} \,\nu(\dd z)\right|\leq&~ \int_\bbr |z|\bone_{\{ |z|\in(0,a] \}} \,\nu(\dd z) =  \int_\bbr |z|^\beta |z|^{1-\beta}\bone_{\{ |z|\in(0,a] \}} \,\nu(\dd z)\\
\leq&~ a^{1-\beta}\int_\bbr |z|^\beta\bone_{\{ |z|\leq 1 \}}\,\nu(\dd z).
\end{align*}

\vspace{0.5\baselineskip}\noindent (3) If $\La$ is symmetric, then the left-hand sides of \eqref{infassym} and \eqref{0assym} are identically zero.\halmos 
\epr

\bpr[of Lemma~\ref{rhogint}]
(1) follows from Lemma~\ref{prop2}, (2) holds because $\rho(\bbr_+\times\bbr^d)=1$ implies that $g-\rho*g=(g-g_\infty)-\rho*(g-g_\infty)$, and (3) is simply Young's inequality.\halmos
\epr

\subsection{Proofs for Section~\ref{distprop}}
\bpr[of Proposition~\ref{gcf}] That $m[X]<\infty$ almost surely can be verified in a similar way to Theorem~\ref{localsol}. Also by essentially the same arguments given there, the stochastic Fubini theorem is applicable for $m[X]$ and the result follows from Proposition~\ref{RR}.\halmos
\epr

\bpr[of Proposition~\ref{SRD}] (1) By assumption the function $h\colon \bbr\times\bbr^d\to\bbr$ defined by $h(t,x)=(g-\rho*g)(t,x)\bone_{\bbr_+}(t)$ belongs to $L^1(\bbr\times\bbr^d)$. Therefore, the claim follows from Young's inequality and the observation that, up to a multiplicative constant, $\cov[X(t,x),X(t+\tau,x+\xi)]$ equals the convolution of $h$ with $h^-$ where $h^-(t,x)=h(-t,-x)$.

\vspace{0.5\baselineskip}\noindent (2) That $X$ is well defined and has a finite second moment, follows from $b_1=0$ and $g-\rho*g\in L^2(\bbr_+\times\bbr^d)$. Moreover, 
as $g-\rho*g$ does not change signs, we have
\begin{align*}
	&~\int_0^\infty\int_{\bbr^d} \left| \int_0^\infty \int_{\bbr^d} (g-\rho*g)(s,y)(g-\rho*g)(s+\tau,y+\xi)\,\dd(s,y)\right|\,\dd(\tau,\xi) \\
	&~\qquad = \int_0^\infty\int_{\bbr^d} \int_0^\infty \int_{\bbr^d} |g-\rho*g|(s,y)|g-\rho*g|(s+\tau,y+\xi)\,\dd(s,y) \,\dd(\tau,\xi)\\
	&~\qquad =  \int_0^\infty \int_{\bbr^d} |g-\rho*g|(s,y)\left(\int_0^\infty\int_{\bbr^d} |g-\rho*g|(s+\tau,y+\xi) \,\dd(\tau,\xi)\right)\,\dd(s,y).
\end{align*} 
Since $g-\rho*g\notin L^1(\bbr_+\times\bbr^d)$, the inner integral is infinite for all $(s,y)$, so the whole integral is infinite as well. This shows that $(\tau,\xi)\mapsto |\cov[X(t,x),X(t+\tau,x+\xi)]|$ is not an element of $L^1(\bbr_+\times\bbr^d)$.
\halmos
\epr

\subsection{Proofs for Section~\ref{pathprop}}\label{proofspocadlag}

\begin{Proof}[of Theorem~\ref{contversion}]
	For $\tilde g:= g-\rho*g$ we have that
	\[ \int_0^t\int_{\bbr^d}\tilde g(t-s,x-y)\,\dd(s,y)=\int_0^t\int_{\bbr^d}\tilde g(s,y)\,\dd(s,y) \]
	is continuous in $(t,x)$. Hence, we may assume without loss of generality that $b=0$.
	The additional conditions on $g$ in Theorem~\ref{contversion}, together with Hölder's inequality and Fubini's theorem, imply for $\rho*g$ (we extend $g$ on the negative half space $(-\infty,0)\times\bbr^d$ by zero):
	\begin{align}
	\notag&\int_0^T\int_{\bbr^d}\left| (\rho*g)(s,y)- (\rho*g)(s+\tau,y+\xi)\right|^2\,\dd(s,y)\\
	\notag&\qquad=\int_0^T\int_{\bbr^d}\bigg| \int_0^s\int_{\bbr^d}  g(s-r,y-z)\,\rho(\dd r,\dd z)\\
	\notag&\qquad \quad-\int_0^{s+\tau}\int_{\bbr^d}  g(s+\tau-r,y+\xi-z)\,\rho(\dd r,\dd z) \bigg|^2\,\dd(s,y)\\
	\notag&\qquad \leq C_T\int_0^T\int_{\bbr^d} \int_0^T\int_{\bbr^d}  \left|g(s-r,y-z)- g(s+\tau-r,y+\xi-z)\right|^2\,\dd(s,y)\,|\rho|(\dd r,\dd z) \\
	\notag&\qquad \quad +C_T\int_0^T\int_{\bbr^d} \int_0^T\int_{\bbr^d}  \left|\bone_{[s,s+\tau]}(r)g(s+\tau-r,y+\xi-z)\right|^2\,\dd(s,y)\,|\rho|(\dd r,\dd z) \\
	\notag&\qquad \leq C_T\int_0^T\int_{\bbr^d} |(\tau,\xi)|^u\,|\rho|(\dd r,\dd z)+C_T\int_0^\tau\int_{\bbr^d} |g(s,y)|^2\,\dd(s,y) \\	
	\notag&\qquad \leq C_T |(\tau,\xi)|^u+C_T\left(\int_0^{\tau}1\,\dd s\right)^{\frac{p-1}{p}}\left(\int_0^{\tau}\left(\int_{\bbr^d} g(s,y)^2\,\dd y\right)^p\,\dd s\right)^\frac{1}{p} \\
	\label{help6}&\qquad \leq C_T |(\tau,\xi)|^u+ C_{T}|(\tau,\xi)|^\frac{p-1}{p},
	\end{align}
	where $|(\tau,\xi)|$ is small enough. Furthermore, we have by another application of Hölder's inequality and Fubini's theorem
	\begin{align}
		\notag\int_0^T\left(\int_{\bbr^d} (\rho*g)(s,y)^2\,\dd y\right)^p\,\dd s & \leq C_T \int_0^T\left(\int_{\bbr^d} \int_0^s\int_{\bbr^d}g(s-r,y-z)^2\,|\rho|(\dd r,\dd z)\,\dd y\right)^p\,\dd s\\
		\notag& \leq C_T \int_0^T\int_0^s\int_{\bbr^d}\left(\int_{\bbr^d} g(s-r,y)^2\,\dd y\right)^p\,|\rho|(\dd r,\dd z)\,\dd s\\
		\label{help7}& \leq C_T \int_0^T\int_{\bbr^d}\int_0^T\left(\int_{\bbr^d} g(s,y)^2\,\dd y\right)^p\,\dd s\,|\rho|(\dd r,\dd z)<\infty.
	\end{align}
	Next, by Corollary~\ref{acf} we have for all $(t,x),(\tau,\xi)\in\bbr_+\times\bbr^d$ that
	\begin{align*}
		\notag&~\bbe\Big[|X(t,x)-X(t+\tau,x+\xi)|^2\Big]\\
		\notag&~\qquad=\si^2\int_0^t\int_{\bbr^d}\tilde g^2(s,y)\,\dd(s,y)-2\si^2\int_0^t\int_{\bbr^d}\tilde g(s,y)\tilde g(s+\tau,y+\xi)\,\dd(s,y)\\
		\notag&~\qquad\quad+\si^2\int_0^{t+\tau}\int_{\bbr^d}\tilde g^2(s,y)\,\dd(s,y)\\
		\notag&~\qquad=2\si^2\int_0^t\int_{\bbr^d}\tilde g(s,y)[\tilde g(s,y)-\tilde g(s+\tau,y+\xi)]\,\dd(s,y)+\si^2\int_t^{t+\tau}\int_{\bbr^d}\tilde g^2(s,y)\,\dd(s,y).
	\end{align*}
	The assumptions on $g$ and the inequalities \eqref{help6} and \eqref{help7} yield
	\[ \int_0^t\int_{\bbr^d}\tilde g(s,y)[\tilde g(s,y)-\tilde g(s+\tau,y+\xi)]\,\dd(s,y) \leq C_T |(\tau,\xi)|^{u/2}+ C_{T}|(\tau,\xi)|^\frac{p-1}{2p} \]
	for small $|(\tau,\xi)|$ and
	\[ \int_t^{t+\tau}\int_{\bbr^d}\tilde g^2(s,y)\,\dd(s,y)\leq\left(\int_t^{t+\tau}1\,\dd s\right)^{\frac{p-1}{p}}\left(\int_t^{t+\tau}\left(\int_{\bbr^d}\tilde g(s,y)^2\,\dd y\right)^p\,\dd s\right)^\frac{1}{p}\leq C_{T}|(\tau,\xi)|^\frac{p-1}{p},  \]
	where we have used the Cauchy-Schwarz inequality in the first and Hölder's inequality in the second step. Adding both inequalities together and using the fact that $X$ is a Gaussian process, Kolmogorov's continuity theorem (see e.g. Theorem~3.23 of \cite{Kallenberg02}) finishes the proof.
	\halmos
\end{Proof}

\begin{Proof}[of Theorem~\ref{tempcadlag}]
	We prove the case when space is one-dimensional, i.e. $d=1$. For higher dimensions the proof is similar.
	Clearly, it suffices to show the path property separately for the drift and Gaussian part, the compensated small jumps part, and the large jumps part.
	
	\vspace{\baselineskip}\noindent
	\textbf{Case 1: $\La(\dd t,\dd x)=b\,\dd (t, x)+\si W(\dd t,\dd x).$}
	
	\vspace{0.5\baselineskip}\noindent	
	The assumptions on $g$ imply that $g\vert_{[0,T]\times\bbr^d}$ belongs to the Sobolev space $W^{1,2}([0,T]\times\bbr^d)$ for all $T\in\bbr_+$ if $\si\neq0$. Therefore, Theorem~3 in Section~5.8 of \cite{Evans98}, a characterization of the Sobolev space $W^{1,2}([0,T]\times\bbr^d)$, ensures the first condition in Theorem~\ref{contversion}. Moreover, since $g,\partial_t g \in L_\loc^2(\bbr_+\times\bbr)$, the fundamental theorem of calculus yields that $s\mapsto\int_{\bbr^d}g^2(s,y)\,\dd y$ is continuous in $s$. Hence, by Theorem~\ref{contversion}, $X$ has a continuous version.
	
	\vspace{\baselineskip}\noindent
	\textbf{Case 2: $\La(\dd t, \dd x )=\int_\bbr z\bone_{\{|z|\leq 1 \}} \,(\pf-\qf)(\dd t,\dd x,\dd z).$}
	
	\vspace{0.5\baselineskip}\noindent	
	We define $\La_n(\dd t, \dd x ):=\int_\bbr z\bone_{\{1/n\leq|z|\leq 1 \}}\bone_{\{|x|\leq n \}} \,(\pf-\qf)(\dd t,\dd x,\dd z)$. Since $\La_n$ has only finitely many jumps on $[0,T]\times\bbr$ almost surely, it is easy to see that the paths of $F_n:=g*\La_n$ are almost surely t-càdlàg and in $L_\loc^\infty(\bbr_+\times\bbr^d)$ due to the boundedness of $g$. Hence, also the process $\rho*F_n$ has almost surely t-càdlàg realizations by dominated convergence. We now show that $\rho*F^n:=\rho*F-\rho*F_n$ converges uniformly on compacts in probability to $0$. To this end, consider for any $T>0$ and $U=[-K,K]\subseteq\bbr$
	\begin{align}
		\notag\bbe\left[\sup_{(t,x)\in[0,T]\times U}\left|(\rho*F^n)(t,x)\right|\right]
		=&~\bbe\left[\sup_{(t,x)\in[0,T]\times U}\left|\int_0^t\int_{\bbr}F^n(t-s,x-y)\,\rho(\dd s,\dd y)\right|\right]\\
		\notag\leq&~\int_0^T\int_{\bbr}\bbe\left[\sup_{(t,x)\in[0,T]\times U}|F^n(t-s,x-y)|\right]\,|\rho|(\dd s,\dd y)\\
		\label{help10}=&~\int_0^T\int_{\bbr}\bbe\left[\sup_{(t,x)\in[0,T]\times U}|F^n(t,x-y)|\right]\,|\rho|(\dd s,\dd y).
	\end{align}
	Now use the fundamental theorem of calculus to decompose $g$ as
	\begin{align}
	\notag	g(t-s,x-\eta-y)=&~g(0,x-\eta-y)+\int_s^t g_1(r-s,x-\eta-y)\,\dd r\\
	\notag	=&~g(0,-\eta-y)+\int_{-\eta}^{x-\eta} g_2 (0,z-y)\,\dd z+\int_s^t g_1(r-s,-\eta-y)\,\dd r \\
	\label{help5}	&~+\int_s^t \int_{-\eta}^{x-\eta} g_{12}(r-s,z-y)\,\dd z\,\dd r,
	\end{align}
	where $\eta\in\bbr$, $g_1(t,x)={\pd_t g(t,x)}$, $g_2(t,x)={\pd_x g(t,x)}$ and $g_{12}(t,x)={\pd_t\pd_x g(t,x)}$. With the same reasoning as in the proof of Theorem~\ref{localsol}, the assumptions imply that the stochastic Fubini theorem is applicable and this gives us
	\begin{align}
	\notag F^n(t,x-\eta):=&~\int_{0}^t\int_{\bbr}g(t-s,x-\eta-y)\,\La^n(\dd s,\dd y)\\
	\notag =&~\int_{0}^t\int_{\bbr}g(0,-\eta-y)\,\La^n(\dd s,\dd y)+\int_{0}^t\int_{\bbr}\int_{-\eta}^{x-\eta} g_2 (0,z-y)\,\dd z\,\La^n(\dd s,\dd y)\\
	\notag &~+\int_{0}^t\int_{\bbr}\int_s^t g_1(r-s,-\eta-y)\,\dd r\,\La^n(\dd s,\dd y)\\ 
	\notag &~+\int_{0}^t\int_{\bbr}\int_s^t \int_{-\eta}^{x-\eta} g_{12}(r-s,z-y)\,\dd z\,\dd r\,\La^n(\dd s,\dd y)\\
	\notag =&~\int_{0}^t\int_{\bbr}g(0,-\eta-y)\,\La^n(\dd s,\dd y)+\int_{-\eta}^{x-\eta}\int_{0}^t\int_{\bbr} g_2 (0,z-y)\,\La^n(\dd s,\dd y)\,\dd z\\
	\notag &~+\int_{0}^t\int_0^r\int_{\bbr} g_1(r-s,-\eta-y)\,\La^n(\dd s,\dd y)\,\dd r\\
	\notag &~+\int_{0}^t\int_{-\eta}^{x-\eta}\int_0^r \int_{\bbr} g_{12}(r-s,z-y)\,\La^n(\dd s,\dd y)\,\dd z\,\dd r\\
	\notag =&\!\!:~I_{1,n}(t,x,\eta)+I_{2,n}(t,x,\eta)+I_{3,n}(t,x,\eta)+I_{4,n}(t,x,\eta),
	\end{align}
	where $\La^n:=\La-\La_n$. Therefore, we have for fixed $\eta\in\bbr$	
	\begin{align*}
		\bbe\left[\sup_{(t,x)\in[0,T]\times U}|F^n(t,x-\eta)|\right]\leq&~\sum_{j=1}^4 \bbe\left[\sup_{(t,x)\in[0,T]\times U}|I_{j,n}(t,x,\eta)|\right].
		% +\bbe\left[\sup_{(t,x)\in[0,T]\times U}|I_{2,n}(t,x)|\right]\\
		%&~+\bbe\left[\sup_{(t,x)\in[0,T]\times U}|I_{3,n}(t,x)|\right]+\bbe\left[\sup_{(t,x)\in[0,T]\times U}|I_{4,n}(t,x)|\right],
	\end{align*}
	Since $I_{1,n}(t,x)$ does not depend on $x$ and is a martingale in $t$, we have by the Burkholder-Davis-Gundy inequalities that	
	\begin{align}
		\notag&\bbe\left[\sup_{(t,x)\in[0,T]\times U}|I_{1,n}(t,x,\eta)|^2\right]\\
		\notag&\qquad=\bbe\Bigg[\sup_{t\in[0,T]}\bigg|\int_0^t\int_\bbr\int_\bbr g(0,-\eta-y)z\Big(\bone_{\{|y|>n\}}\bone_{\{|z|\leq1\}}\\
		\notag&\qquad\quad+\bone_{\{|y|\leq n\}}\bone_{\{|z|<1/n\}}\Big)(\pf-\qf)(\dd s,\dd y,\dd z)\bigg|^2\Bigg] \\
		\notag&\qquad\leq C_T\bbe\left[\int_0^T\int_\bbr\int_\bbr |g|^2(0,-\eta-y)|z|^2\left(\bone_{\{|y|>n\}}\bone_{\{|z|\leq1\}}+\bone_{\{|y|\leq n\}}\bone_{\{|z|<1/n\}}\right)\pf(\dd s,\dd y,\dd z)\right] \\
		\notag&\qquad=C_T\int_0^T\int_\bbr |g|^2(0,-\eta-y)\bone_{\{|y|>n\}}\,\dd(s,y)\int_\bbr|z|^2\bone_{\{|z|\leq1\}}\,\nu(\dd z)\\
		\label{help8}&\qquad\quad+C_T\int_0^T\int_\bbr |g|^2(0,-\eta-y)\bone_{\{|y|\leq n\}}\,\dd(s,y)\int_\bbr|z|^2\bone_{\{|z|<1/n\}}\,\nu(\dd z)\to0
	\end{align} 
	as $n\to\infty$ and that $\bbe\left[\sup_{(t,x)\in[0,T]\times U}|I_{1,n}(t,x,\eta)|^2\right]$ is bounded in $n$ and $\eta$.
	Here we have used that $g\in L^2_\loc(\bbr_+\times\bbr)$ since it is continuous on $\bbr_+\times\bbr$ and belongs to $L^\al_\loc(\bbr_+\times\bbr)$. 
	By similar arguments, we have that
	\begin{align*}
		&\bbe\left[\sup_{(t,x)\in[0,T]\times U}|I_{2,n}(t,x,\eta)|^2\right]\\
		&\qquad=\bbe\left[\sup_{(t,x)\in[0,T]\times U}\left|\int_{-\eta}^{x-\eta}\int_{0}^t\int_{\bbr} g_2 (0,z-y)\,\La^n(\dd s,\dd y)\,\dd z\right|^2\right]\\
		&\qquad\leq C_T\int_{-K-\eta}^{K-\eta}\bbe\left[\sup_{(t,x)\in[0,T]\times U}\left|\int_{0}^t\int_{\bbr} g_2 (0,z-y)\,\La^n(\dd s,\dd y)\right|^2\right]\,\dd z\\
		&\qquad\leq C_T\int_{-K-\eta}^{K-\eta}\int_0^T\int_\bbr |g_2|^2(0,z-y)\bone_{\{|y|>n\}}\,\dd(s,y)\int_\bbr|\zeta|^2\bone_{\{|\zeta|\leq1\}}\,\nu(\dd \zeta)\,\dd z\\
		&\qquad\quad+C_T\int_{-K-\eta}^{K-\eta}\int_0^T\int_\bbr |g_2|^2(0,z-y)\bone_{\{|y|\leq n\}}\,\dd(s,y)\int_\bbr|\zeta|^2\bone_{\{|\zeta|<1/n\}}\,\nu(\dd \zeta)\,\dd z\\
		&\qquad\leq C_T\int_0^T\int_\bbr |g_2|^2(0,y)\bone_{\{|y|>n-K-|\eta|\}}\,\dd(s,y)\int_\bbr|\zeta|^2\bone_{\{|\zeta|\leq1\}}\,\nu(\dd \zeta)\\
		&\qquad\quad+C_T\int_0^T\int_\bbr |g_2|^2(0,y)\,\dd(s,y)\int_\bbr|\zeta|^2\bone_{\{|\zeta|<1/n\}}\,\nu(\dd \zeta)\to 0
	\end{align*}
	as $n\to 0$ and that $\bbe\left[\sup_{(t,x)\in[0,T]\times U}|I_{2,n}(t,x,\eta)|^2\right]$ is bounded in $n$ and $\eta$.
	Because $I_{3,n}$ and $I_{4,n}$ can be treated analogously to $I_{2,n}$, $\rho*F_n$ converges uniformly on compacts in probability to $\rho*F$ due to \eqref{help10} and dominated convergence. This gives us a t-càdlàg version of $\rho*F$. By setting $\eta=0$, we also obtain that $F_n$ converges uniformly on compacts in probability to $F$. As a consequence, $X$ has a t-càdlàg version.
	
	\vspace{0.5\baselineskip}\noindent
	\textbf{Case 3: $\La(\dd t, \dd x )=\int_\bbr z\bone_{\{|z|> 1 \}} \,\pf(\dd t,\dd x,\dd z).$} 
	
	\vspace{0.5\baselineskip}\noindent	
	Here we assume $\al<1$ because in the situation $\al=1$ we can split the Lévy basis according to 
	\[  \La(\dd t, \dd x )=\int_\bbr z\bone_{\{|z|> 1 \}} \,(\pf-\qf)(\dd t,\dd x,\dd z)+\int_\bbr z\bone_{\{|z|> 1 \}} \,\qf(\dd t,\dd x,\dd z),\]
	treating the first summand as in Case 2 and the second summand as in Case 1. We consider the truncated Lévy basis $\La_n(\dd t, \dd x ):=\int_\bbr z\bone_{\{|z|\geq 1 \}}\bone_{\{|x|\leq n \}} \,\pf(\dd t,\dd x,\dd z)$, which almost surely has finitely many jumps on $[0,T]\times\bbr$. %Since $X=F-\rho*F$ (where $F=g*\La$), it is enough to show that both $F$ and $\rho*F$ have a t-càdlàg version. 
	As in Case 2, the processes $F_n=g*\La_n$ and $\rho*F_n$ have t-càdlàg paths almost surely. It suffices therefore to prove that they converge uniformly on compacts in probability to $F$ and $\rho*F$, respectively. 
	We can estimate (note that we can interchange convolution with $\La$ and convolution with $\rho$ as shown in Theorem~\ref{localsol})
	%\[ \rho*F_n(t,x)=\int_0^t\int_\bbr\int_\bbr(\rho*g)(t-s,x-y)z \bone_{\{|z|\geq 1 \}}\bone_{\{|y|\leq n \}} \,\pf(\dd s,\dd y,\dd z),\]	
	\begin{align*}
		&\bbe\left[\sup_{(t,x)\in[0,T]\times U}|(\rho*F-\rho*F_n)(t,x)|^\al\right]\\
		&\qquad = \bbe\left[\sup_{(t,x)\in[0,T]\times U}\left|\int_0^t\int_{\bbr} (\rho*g)(t-s,x-y)\,(\La-\La_n)(\dd s,\dd y)\right|^\al\right]\\
		&\qquad\leq \int_0^T\int_{\bbr}\sup_{(t,x)\in[0,T]\times U} \left|(\rho*g)\right|^\al(t-s,x-y)\bone_{\{|y|> n \}}\,\dd (s,y)\int_\bbr |z|^\al\bone_{\{|z|\geq 1 \}}\,\nu(\dd z)\\
		&\qquad\leq C_T \int_{\bbr}\sup_{(t,x)\in[0,T]\times U} \left|(\rho*g)\right|^\al(t,x-y)\bone_{\{|y|> n \}}\,\dd y\\
		&\qquad= C_T\int_{\bbr}\sup_{(t,x)\in[0,T]\times U} \left|\vp(-y)(\rho*g)(t,x-y)\right|^\al\bone_{\{|y|> n \}}\vp(-y)^{-\al}\,\dd y,		
	\end{align*}
	where $\vp$ is the function from Theorem~\ref{localsol}. Since $\vp^{-\al}\in L^1(\bbr)$, the right-hand side of the last inequality tends to zero by dominated convergence if we can show that 
	\[ \sup_{(t,x)\in[0,T]\times U} \vp(-y)|\rho*g|(t,x-y) \]
	is bounded in $y$. But this follows because for every $y\in\bbr$ we have
	\begin{align*}
		\sup_{(t,x)\in[0,T]\times U} \vp(-y)|\rho*g|(t,x-y)\leq&~\sup_{(t,x)\in[0,T]\times U} \vp(x-y)|\rho*g|(t,x-y)\vp(-x)\\
		\leq&~C_T\sup_{x\in U}\vp(-x)<\infty,
	\end{align*}
	where we have used that $\vp(|\rho|*|g|) \in L^\infty_\loc(\bbr_+\times\bbr)$ as shown in the proof of Theorem~\ref{localsol} and that $\vp$ is submultiplicative and locally bounded. %Since $\rho*F_n$ has t-càdlàg paths almost surely, we can choose a t-càdlàg version for $\rho*F$. 
	Similar arguments applied to the pair $F_n$ and $F$ yield
	\begin{align}
	\notag\bbe\left[\sup_{(t,x)\in[0,T]\times U}|(F-F_n)(t,x)|^\al\right]
	&\leq C_T\int_{\bbr}\sup_{(t,x)\in[0,T]\times U} \left|g\right|^\al(t,x-y)\bone_{\{|y|> n \}}\,\dd y\\
	\notag&\leq C_T\int_{\bbr}\sup_{x\in U} G^\al(|x-y|)\bone_{\{|y|> n \}}\,\dd y\\
	\label{help9}&\leq C_T\int_{\bbr}G^\al((|y|-K)\vee0)\bone_{\{|y|> n \}}\,\dd y,
	\end{align}
	where we used the monotonicity of $G$. Now the last line goes to $0$ by dominated convergence because $G(|x|)\in L^\al(\bbr)$. Altogether we obtain a t-càdlàg version of $X$.
	\halmos
\end{Proof}

For the proof of Theorem~\ref{pocadlag} we %use another strategy: the desired path property is first proved for a sequence of approximating processes. Subsequently we show that this sequence converges uniformly on compact sets to our solution process. As a consequence the path property is transferred to $X$. In order to establish uniform convergence we are going to 
need to resort to maximal inequalities for multi-parameter martingales in line with \cite{Khoshnevisan02,Walsh79}. 
To this end, let $\leq$ denote the partial order on $\bbr^2$ such that $v=(v_1,v_2)\leq w=(w_1,w_2)$ if and only if $v_1\leq w_1$ and $v_2\leq w_2$. For two subsets $I_1$ and $I_2$ of $\bbr$ set $I=I_1\times I_2$. Now a family of sub-$\si$-algebras $\bbg=(\calg(v))_{v\in I}$ is called a filtration if $\calg(v)\subseteq\calg(w)$ for all $v\leq w$ in $I$. A stochastic process $X$ indexed by $I$ is called a \emph{martingale} with respect $\bbg$ if $X$ is adapted to $\bbg$, $X(v)$ is integrable for all $v\in I$ and $\bbe[X(w)\mid\calg(v)]=X(v)$ for all $v\leq w$ in $I$. Furthermore, we define the marginal filtrations $\bbg^1=(\calg^1(v_1))_{v_1\in I_1}$ and $\bbg^2=(\calg^2(v_2))_{v_2\in I_2}$ through $\calg^1(v_1):=\bigvee_{\xi\in I_2}\calg(v_1,\xi)$ and $\calg^2(v_2):=\bigvee_{\xi\in I_1}\calg(\xi,v_2)$ and set $\calg^*(v_1,v_2):=\calg^1(v_1)\vee\calg^2(v_2)$. Then a martingale $X$ is called an \emph{orthomartingale} if for each $(i,j)\in\{(1,2),(2,1)\}$ and each fixed $v_i\in I_i$, $v_j\mapsto X(v)$ is a one-parameter martingale with respect to $\bbg^j$. Moreover, a martingale $X$ is called a \emph{strong martingale} if it satisfies the condition $\bbe[X((v_1,v_2),(w_1,w_2)]\mid\calg^*(v_1,v_2)]=0$ for all $(v_1,v_2)\leq(w_1,w_2)$ in $I$, where $X((v_1,v_2),(w_1,w_2)]:=X(w_1,w_2)-X(w_1,v_2)-X(v_1,w_2)+X(v_1,v_2)$ is the two-dimensional increment. Further notation includes $[v,w]_\leq:=\{u\in\bbr^2\colon v\leq u\leq w\}$ for the closed interval from $v$ to $w$ with respect to the partial order $\leq$. Similarly, for $(t,x),(\tilde t,\tilde x)\in\bbr_+\times\bbr$, $[(t,x),(\tilde t,\tilde x)]_\lpo:=\{(s,y)\in\bbr_+\times\bbr\colon(t,x)\lpo(s,y)\lpo(\tilde t,\tilde x)\}$ denotes the closed interval from $(t,x)$ to $(\tilde t,\tilde x)$ with respect to the partial order $\lpo$ as defined in Section~\ref{pathprop}. Also, we use the abbreviations $A:=-\{(t,x)\in\bbr_+\times\bbr\colon|x|\leq ct\}$, $A(t,x):=A+(t,x)$ and $A^+(t,x):=A(t,x)\cap(\bbr_+\times\bbr)$. 

The following lemma extends the previously known maximal inequalities for multi-parameter martingales \cite{Khoshnevisan02, Walsh79} to processes that are not martingales themselves but can be seen as ``rotated martingales''. For later purposes, we also need the situation where $\La$ is a not necessarily homogeneous Lévy basis (i.e., the coefficients $b$, $\si$ and $\nu$ in \eqref{levybasis} may depend on $(t,x)$ in such a way that $\La(A)$ is well defined for all $A\in\calb_\bb(\bbr_+\times\bbr^d)$).

\begin{Lemma}\label{maxinequalities}
	\begin{enumerate}
		\item If $\La$ is a (not necessarily homogeneous) Lévy basis with mean $0$, the process 
		\[ X(t,x)=\int_{0}^t\int_{\bbr}\bone_{A(t,x)}(s,y)\,\La(\dd s,\dd y)=\La(A^+(t,x)),\quad (t,x)\in\bbr_+\times\bbr, \]
		satisfies the maximal inequality
		\begin{equation*}
		\la\bbp\left[\sup_{(s,y)\in[(\tilde t,\tilde x),(t,x)]_\lpo}|X(s,y)|\geq\la\right]\leq13\bbe\left[|X(t,x)|\right]
		\end{equation*}
		for every $\la>0$ and $(\tilde t,\tilde x)\lpo(t, x)$ in $\bbr_+\times\bbr$.
		\item If $\La$ further has a finite $p$'th moment with some $p>1$, then $X$ satisfies 
		\begin{equation*}
		\bbe\left[\sup_{(s,y)\in[(\tilde t,\tilde x),(t,x)]_\lpo}|X(s,y)|^p\right]\leq\left(\frac{p}{p-1}\right)^{2p}\bbe[|X(t,x)|^p]
		\end{equation*}
		for every $(\tilde t,\tilde x)\lpo(t, x)$ in $\bbr_+\times\bbr$.		
	\end{enumerate}
\end{Lemma}
\begin{Proof}
	Without loss of generality we assume $(\tilde t,\tilde x)=(0,0)$ and define a family of sub-$\si$-algebras $\bbf=(\calf(s,y))_{(s,y)\in [(0,0),(t,x)]_\lpo}$ by
	\[ \calf(s,y):=\sigma(X(\tau,\xi)\colon(\tau,\xi)\in\bbr_+\times\bbr,(\tau,\xi)\lpo(s, y)). \]
	Then we have $\calf(\tilde s,\tilde y)\subseteq\calf(s, y)$ for every $(\tilde s,\tilde y)\lpo(s, y)$ and the process $X$ is integrable and adapted to $\bbf$ on $[(0,0),(t,x)]_\lpo$. Next, we define for any $(\tilde s,\tilde y)\lpo(s, y)$ in $[(0,0),(t,x)]_\lpo$
	\[ \calm(\tilde s, \tilde y):=\big\{\{X(s_1,y_1)\in B_1,...,X(s_n,y_n)\in B_n\}\colon n\in\bbn,(s_i,y_i)\lpo(\tilde s,\tilde y),B_i\in\calb(\bbr)\big\}.\]
	Then the properties of a Lévy basis imply $X(s,y)-X(\tilde s,\tilde y)\indep\calm(\tilde s, \tilde y)$. 
	Since $\calm(\tilde s, \tilde y)$ is intersection-stable, we get $X(s,y)-X(\tilde s,\tilde y)\indep\sigma(\calm(\tilde s, \tilde y))=\calf(\tilde s,\tilde y)$ and therefore
	\begin{align}
		\notag\bbe[X(s,y)\mid\calf(\tilde s,\tilde y)]&=\bbe[\La(A^+(\tilde s,\tilde y))+\La(A^+(s,y)\backslash A^+(\tilde s,\tilde y))\mid\calf(\tilde s,\tilde y)]\\
		\notag&=\La(A^+(\tilde s,\tilde y))+\bbe[\La(A^+(s,y)\backslash A^+(\tilde s,\tilde y))\mid\calf(\tilde s,\tilde y)]\\
		\label{help1}&=\La(A^+(\tilde s,\tilde y))=X(\tilde s,\tilde y).
	\end{align}
	We now transform $X$ into a strong martingale by considering the function $H\colon[(0,0),(t,x)]_\lpo\to [(0,0),H(t,x)]_\leq$  given by
	\[ H(s,y)=\begin{pmatrix}
	\frac{1}{\sqrt{2}}&-\frac{1}{\sqrt{2}}\\
	\frac{1}{\sqrt{2}}&\frac{1}{\sqrt{2}}
	\end{pmatrix}
	\begin{pmatrix}1&0\\0&\frac{1}{c}\end{pmatrix}
	\begin{pmatrix}s\\y\end{pmatrix}. \]
	Note that the first matrix is the rotation matrix about 45 degrees counter-clockwise and the second matrix is a rescaling in the space coordinate. It is easy to see that $H$ is in fact order-preserving and bijective. Now let $\tilde X$ be the push-forward process of $X$ through $H$, i.e. $\tilde X(v_1,v_2):=X(H^{-1}(v_1,v_2))$, and $\tilde \bbf$ be the push-forward of $\bbf$ through $H$, i.e. $\tilde \calf(v_1,v_2):=\calf(H^{-1}(v_1,v_2))$ for every $(v_1,v_2)\in[(0,0),H(t,x)]_\leq$. Then $\tilde \bbf$ is a filtration on $[(0,0),H(t,x)]_\leq$ and  $\tilde X$ is a martingale with respect to $\tilde\bbf$ since the property in \eqref{help1} of $X$ is inherited through $H$. In fact, $\tilde X$ is even a strong martingale because we have for $(v_1,v_2)\leq(w_1,w_2)$ in $[(0,0),H(t,x)]_\leq$ that
	\begin{align*}
		\tilde X((v_1,v_2),(w_1,w_2)]&=\tilde X(w_1,w_2)-\tilde X(v_1,w_2)-\tilde X(w_1,v_2)+\tilde X(v_1,v_2)\\
		& =X(H^{-1}(w_1,w_2))-X(H^{-1}(v_1,w_2))-X(H^{-1}(w_1,v_2))+X(H^{-1}(w_1,w_2))\\
		& = \La(A^+(H^{-1}(w_1,w_2)))-\La(A^+(H^{-1}(v_1,w_2)))-\La(A^+(H^{-1}(w_1,v_2)))\\
		& \quad +\La(A^+(H^{-1}(w_1,w_2)))\\
		&=\La([H^{-1}(v_1,v_2),H^{-1}(w_1,w_2)]_\lpo),	
	\end{align*} 
	where the last inequality follows from the triangular shape of $A^+$. Moreover, letting $(u_1,u_2):=H(t,x)$, we have
	\begin{align*}
		\tilde\calf^*(v_1,v_2)&=\tilde\calf^1(v_1)\vee\tilde\calf^2(v_2)=\tilde\calf(v_1,u_2)\vee\tilde\calf(u_1,v_2)=\calf(H^{-1}(v_1,u_2))\vee\calf(H^{-1}(u_1,v_2))\\
		&=\sigma\left(X(s,y)\colon (s,y)\in A^+(H^{-1}(v_1,u_2))\cup A^+(H^{-1}(u_1,v_2))\right).
	\end{align*}
	With the same argument as in \eqref{help1}, we can show that $\La([H^{-1}(v_1,v_2),H^{-1}(w_1,w_2)]_\lpo)$ is independent of 
	$\sigma(X(s,y)\colon (s,y)\in A^+(H^{-1}(v_1,u_2))\cup A^+(H^{-1}(u_1,v_2)))$, which implies	
	\begin{equation*}
		\bbe[\tilde X((v_1,v_2),(w_1,w_2)]\mid\tilde\calf^*(v_1,v_2)]=0.
	\end{equation*}	
	Therefore $\tilde X$ is a strong martingale with respect to $\tilde \bbf$. This allows us to use Walsh's maximal inequality for strong martingales: by Corollary 3.4 in \cite{Walsh79} we get
	\begin{equation*}
		\la\bbp\left[\sup_{(s,y)\in[(0,0),H(t,x)]_\leq}|\tilde X(u_1,u_2)|\geq\la\right]\leq13\sup_{(s,y)\in[(0,0),H(t,x)]_\leq}\bbe[|\tilde X(u_1,u_2)|]=13\bbe[|\tilde X(H(t,x))|]
	\end{equation*}	
	for all $\la>0$. By the definition of $\tilde X$ this is equivalent to
	\begin{equation*}
		\la\bbp\left[\sup_{(s,y)\in[(0,0),(t,x)]_\lpo}|X(s,y)|\geq\la\right]\leq13\bbe[|X(t,x)|].
	\end{equation*}
	The second part of the lemma can be proved similarly, using Cairoli's maximal inequality for orthomartingales (see e.g. Corollary 2.3.1 in \cite{Khoshnevisan02}) and the fact that a strong martingale is always an orthomartingale (see e.g. Proposition 1.1 in \cite{Walsh79}).
	\halmos
\end{Proof}

\begin{Proof}[of Theorem~\ref{pocadlag}]
	
	The cases where $\La$ is equal to the drift and Gaussian part, the compensated small jumps part, or the large jumps part are considered separately.
	
	\vspace{\baselineskip}\noindent
	\textbf{Case 1: $\La(\dd t,\dd x)=b\,\dd (t, x)+\si W(\dd t,\dd x).$}
	
	\vspace{0.5\baselineskip}\noindent
	Our assertion is proved once we can show both conditions of Theorem~\ref{contversion}. The second condition is obviously satisfied and the first condition follows similarly as in Case 1 of the proof of Theorem~\ref{tempcadlag} in conjunction with the boundedness of $h$. We omit the details here.
	
	\vspace{\baselineskip}\noindent
	\textbf{Case 2: $\La(\dd t, \dd x )=\int_\bbr z\bone_{\{|z|\leq 1 \}} \,(\pf-\qf)(\dd t,\dd x,\dd z).$}
	
	\vspace{0.5\baselineskip}\noindent
	The argument is similar to Case 2 of the proof of Theorem~\ref{tempcadlag}. Therefore, we only highlight the major differences.	
	\begin{comment}
	gives us the following decomposition of $h$:
	\begin{align*}
		h(t-s,x-y)=&~h(0,x-y)+\int_s^t h_1(r-s,x-y)\,\dd r\\
		=&~h(0,-y)+\int_0^x h_2 (0,z-y)\,\dd z\\
		&~+\int_s^t h_1(r-s,-y)\,\dd r +\int_s^t \int_0^x h_{12}(r-s,z-y)\,\dd z\,\dd r,
	\end{align*}
	where $h_1(t,x)={\pd h(t,x)}/{\pd t}$, $h_2(t,x)={\pd h(t,x)}/{\pd x}$ and $h_{12}(t,x)={\pd^2 h(t,x)}/{(\pd t\,\pd x)}$. With the same argument as in the proof of Theorem~\ref{localsol} we conclude that the stochastic Fubini theorem in Theorem~2 of  
	\cite{Lebedev96} is applicable and this gives us
	\end{comment}
	Again, we take advantage of the fact that the Lévy basis $\La_n(\dd t, \dd x ):=\int_\bbr z\bone_{\{1/n\leq|z|\leq 1 \}}\bone_{\{|x|\leq n \}} \,(\pf-\qf)(\dd t,\dd x,\dd z)$ has only finitely many jumps on $[0,T]\times\bbr$ almost surely. Therefore, both $F_n:=g*\La_n$ and $\rho*F_n$ have $\lpo$-càdlàg paths almost surely. Our claim is proved if we show that $F^n:=F-F_n$ and $\rho*F^n=\rho*F-\rho*F_n$ both converge to $0$, uniformly on compacts in probability. To this end, we estimate for $(t_1,x_1)\lpo (t_2,x_2)$ in $\bbr_+\times\bbr$ and a sufficiently big $T>0$
	\small
	\begin{equation*}
		\bbe\left[\sup_{(t,x)\in[(t_1,x_1),(t_2,x_2)]_\lpo}\left|(\rho*F^n)(t,x)\right|\right]
		%=&~\bbe\left[\sup_{(t,x)\in[(t_1,x_1),(t_2,x_2)]_\lpo}\left|\int_0^t\int_{\bbr}F^n(t-s,x-y)\,\rho(\dd s,\dd y)\right|\right]\\
		%\leq&~\int_0^T\int_{\bbr}\bbe\left[\sup_{(t,x)\in[(t_1,x_1),(t_2,x_2)]_\lpo}|F^n(t-s,x-y)|\right]\,|\rho|(\dd s,\dd y)\\
		\leq\int_0^T\int_\bbr\bbe\left[\sup_{(t,x)\in[(t_1,x_1),(t_2,x_2)]_\lpo}|F^n(t-s,x-y)|\right]\,|\rho|(\dd s,\dd y)	
	\end{equation*}
	\normalsize	
	and use the fundamental theorem of calculus and the stochastic Fubini theorem to split the last integrand into four parts according to 
	\begin{align}
	\notag F^n(t-u,x-\eta)=&~\int_{0}^{t-u}\int_\bbr \bone_{A(t-u,x-\eta)}(s,y)h(t-u-s,x-\eta-y)\,\La^n(\dd s,\dd y)\\
	%\notag =&~\int_{0}^t\int_{x-c(t-s)}^{x+c(t-s)}h(0,-y)\,\La(\dd s,\dd y)+\int_{0}^t\int_{x-c(t-s)}^{x+c(t-s)}\int_0^x h_2 (0,z-y)\,\dd z\,\La(\dd s,\dd y)\\
	%\notag &~+\int_{0}^t\int_{x-c(t-s)}^{x+c(t-s)}\int_s^t h_1(r-s,-y)\,\dd r\,\La(\dd s,\dd y)\\ 
	%\notag &~+\int_{0}^t\int_{x-c(t-s)}^{x+c(t-s)}\int_s^t \int_0^x h_{12}(r-s,z-y)\,\dd z\,\dd r\,\La(\dd s,\dd y)\\
	\notag =&~\int_{0}^{t-u}\int_{\bbr} \bone_{A(t-u,x-\eta)}(s,y)h(0,-\eta-y)\,\La^n(\dd s,\dd y)\\
	\notag &~+\int_{-\eta}^{x-\eta}\int_{0}^{t-u}\int_{\bbr} \bone_{A(t-u,x-\eta)}(s,y)h_2 (0,z-y)\,\La^n(\dd s,\dd y)\,\dd z\\
	\notag &~+\int_{0}^{t-u}\int_0^r\int_{\bbr} \bone_{A(t-u,x-\eta)}(s,y)h_1(r-s,-\eta-y)\,\La^n(\dd s,\dd y)\,\dd r \\
	\label{help2} &~+\int_{0}^{t-u}\int_{-\eta}^{x-\eta}\int_0^r \int_{\bbr} \bone_{A(t-u,x-\eta)}(s,y)h_{12}(r-s,z-y)\,\La^n(\dd s,\dd y)\,\dd z\,\dd r,
	\end{align}	
	where $(u,\eta)\in[0,t]\times\bbr$, $h_1(t,x)={\pd_t h(t,x)}$, $h_2(t,x)={\pd_x h(t,x)}$, $h_{12}(t,x)={\pd_t\pd_x h(t,x)}$ and $\La^n=\La-\La_n$. 
	For the first summand, we have by Lemma~\ref{maxinequalities} and a similar reasoning as in \eqref{help8} that
	\begin{align*}
		&\bbe\left[\sup_{(t,x)\in[(t_1,x_1),(t_2,x_2)]_\lpo}\left| \int_{0}^{t-u}\int_\bbr \bone_{A(t-u,x-\eta)}(s,y)h(0,-\eta-y)\,\La^n(\dd s,\dd y) \right|^2\right]\\
		&\qquad\leq C_T\bbe\left[\left| \int_{0}^{t_2-u}\int_\bbr \bone_{A(t_2-u,x_2-\eta)}(s,y)h(0,-\eta-y)\,\La^n(\dd s,\dd y) \right|^2\right]
	\end{align*}
	converges to $0$ as $n\to\infty$ and is bounded in $(u,\eta)$ and $n$. Note that we have used that $h$ belongs to $L_\loc^2(\bbr_+\times\bbr)$ since it is a bounded function in $L_\loc^\al(\bbr_+\times\bbr)$.  Regarding the second summand  we have for sufficiently big $\xi>0$ that
	\begin{align*}
	&\bbe\left[\sup_{(t,x)\in[(t_1,x_1),(t_2,x_2)]_\lpo}\left| \int_{-\eta}^{x-\eta}\int_{0}^{t-u}\int_\bbr \bone_{A(t-u,x-\eta)}(s,y)h_2 (0,z-y)\,\La^n(\dd s,\dd y)\,\dd z \right|^2\right]\\
	&\qquad\leq C_T\int_{-\xi-\eta}^{\xi-\eta}\bbe\left[\sup_{(t,x)\in[(t_1,x_1),(t_2,x_2)]_\lpo}\left| \int_{0}^{t-u}\int_\bbr \bone_{A(t-u,x-\eta)}(s,y)h_2 (0,z-y)\,\La^n(\dd s,\dd y) \right|^2\right]\,\dd z\\
	&\qquad\leq C_T\int_{-\xi-\eta}^{\xi-\eta}\bbe\left[\left| \int_{0}^{t_2-u}\int_\bbr \bone_{A(t_2-u,x_2-\eta)}(s,y)h_2(0,z-y)\,\La^n(\dd s,\dd y) \right|^2\right]\,\dd z\to0,
	\end{align*}
	where we have used Lemma~\ref{maxinequalities} and that $h_2$ also belongs to $L_\loc^2(\bbr_+\times\bbr)$.
	Analogous calculations hold for the third and fourth summand. As a consequence, both $F^n$ and $\rho*F^n$ converge uniformly on compacts in probability to $0$ by dominated convergence.
	
	\vspace{\baselineskip}\noindent
	\textbf{Case 3: $\La(\dd t, \dd x )=\int_\bbr z\bone_{\{|z|> 1 \}} \,\pf(\dd t,\dd x,\dd z).$}

	\vspace{0.5\baselineskip}\noindent
	In this case the same argument as in Case 3 of the proof of Theorem~\ref{tempcadlag} applies with $\bone_{-A}h$ instead of $g$. Notice that, under the current setting, the first integral in \eqref{help9} is actually taken for a bounded function on a compact subset of $\bbr$ (due to the indicator $\bone_{-A}$), is therefore finite and tends to zero as $n\to\infty$. The only difference appears in the case $\al=1$, where we cannot copy the proof of Theorem~\ref{tempcadlag} since Lemma~\ref{maxinequalities}(2) requires $p>1$. Instead, we observe for any $T>0$ and $U=[-K,K]\subseteq\bbr$ that
	\begin{align}
	\notag&\bbe\left[\sup_{(t,x)\in[0,T]\times U}|(\rho*F-\rho*F_n)(t,x)|\right]\\
	\notag&\qquad = \bbe\left[\sup_{(t,x)\in[0,T]\times U}\left|\int_0^t\int_{\bbr} (\rho*g)(t-s,x-y)\,(\La-\La_n)(\dd s,\dd y)\right|\right]\\
	\notag&\qquad\leq \int_0^T\int_{\bbr}\sup_{(t,x)\in[0,T]\times U} \left|(\rho*g)\right|(t-s,x-y)\bone_{\{|y|> n \}}\,\dd (s,y)\int_\bbr |z|\bone_{\{|z|\geq 1 \}}\,\nu(\dd z)\\
	\notag&\qquad\leq C_T \int_{\bbr}\sup_{(t,x)\in[0,T]\times U} \left|(\rho*g)\right|(t,x-y)\bone_{\{|y|> n \}}\,\dd y\\
	\notag&\qquad\leq C_T \int_{\bbr}\sup_{(t,x)\in[0,T]\times U} \int_0^t\int_\bbr |g(t-u,x-y-v)|\,|\rho|(\dd u, \dd v)\bone_{\{|y|> n \}}\,\dd y\\
	\label{help4}&\qquad\leq C_T \int_0^T\int_\bbr \int_{\bbr}\sup_{(t,x)\in[0,T]\times U} |g(t-u,x-y-v)|\bone_{\{|y|> n \}}\,\dd y\,|\rho|(\dd u, \dd v),		
	\end{align}
	and that $\int_{\bbr}\sup_{(t,x)\in[0,T]\times U} |g(t-u,x-y-v)|\bone_{\{|y|> n \}}\,\dd y$ is bounded in $(n,u,v)$ and goes to $0$ by dominated convergence since $g=\bone_{-A}h$. Consequently also \eqref{help4} tends to $0$ by dominated convergence. The pair $F$ and $F_n$ can be treated analogously.
	\halmos	
\end{Proof}	

\begin{appendix}
\section{Proof of Proposition~\ref{thmdetsolution}}\label{AppA}
We first collect some useful properties of convolutions. The proof of the following two lemmata is analogous to the one-parameter case  (see Section~4.1 of \cite{Gripenberg90} or Example~10.3 of \cite{Rudin91} for the first result, and Section~3.6 of \cite{Gripenberg90} for the second result).
\blem\label{prop1}
	Let $S$ be $\bbr^{d+1},\ \bbr_+\times\bbr^d$ or $[0,T]\times\bbr^d$ and $\mu,\ \eta$ and $\pi$ be measures in $M(S)$. Then
	\begin{enumerate}
		\item $\mu*\eta\in M(S)$ and $\|\mu*\eta\|\leq\|\mu\|\|\eta\|$,
		\item $(\mu*\eta)*\pi=\mu*(\eta*\pi)$,
		\item $\mu*\eta=\eta*\mu$.
	\end{enumerate}
The statement is still valid if $M(S)$ is replaced by $M_{\loc}(\bbr_+\times\bbr^d)$ and $\|\cdot\|$ in (1) is replaced by the total variation norm on $[0,T]\times\bbr^d$ for some arbitrary $T\in\bbr_+$.
\elem
\blem\label{prop2}
	Let $\mu$ and $\eta$ be measures in $M_{\loc}(\bbr_+\times\bbr^d)$ and $h\in L_{\loc}^p (\bbr_+\times\bbr^d)$. Then the following statements hold for all $p\in[1,\infty]$.
	\begin{enumerate}
		\item $h*\mu\in L_{\loc}^p (\bbr_+\times\bbr^d)$.
		\item $(h*\mu)*\eta=h*(\mu*\eta)$ and $(\mu*h)*\eta=\mu*(h*\eta)$.
		\item If additionally $\mu\in M(\bbr_+\times\bbr^d)$ and $h\in L^p (\bbr_+\times\bbr^d)$, then $h*\mu$ also belongs to $L^p (\bbr_+\times\bbr^d)$.
	\end{enumerate}
\elem

\bpr[of Proposition~\ref{thmdetsolution}] (1)
	Our proof extends Theorem~4.1.5 of \cite{Gripenberg90}, and for the reader's convenience we present the details in short. Note that this part can alternatively be proven in a more abstract framework using Theorem 4.3.6(b) of \cite{Palmer94} involving the Jacobson radical of the commutative Banachalgebra $(M([0,T]\times\bbr^d),+,*)$. 
	First we show that for each positive $T$ there is a unique $\rho_T$ in $M([0,T]\times\bbr^d)$ such that 
	\[ \rho_T + \mu_T =\mu_T*\rho_T.  \]
	Here $\mu_T$ is the restriction of $\mu$ on $[0,T]\times\bbr^d$.
	To show the existence of $\rho_T$ we construct a geometric series and use a Banach space argument.
	We first consider the special case $\|\mu_T\|=|\mu_T|([0,T]\times\bbr^d)<1$.
	Defining
	\[ \rho_m:=-\sum_{j=1}^{m}\mu_T^{*j},\quad m\in\bbn, \]
	we obtain
	\[ \rho_m+\mu_T=-\sum_{j=1}^{m}\mu_T^{*j}+\mu_T=-\sum_{j=2}^{m}\mu_T^{*j}=\mu_T*\left(-\sum_{j=1}^{m-1}\mu_T^{*j}\right)=\mu_T*\rho_{m-1},\quad m\in\bbn\backslash\{1\}. \]
	By Lemma~\ref{prop1}, we have $\|\mu_T^{*j}\|\leq\|\mu_T\|^j$,
	so $(\rho_m)$ is a Cauchy sequence and converges to some $\rho_T\in M([0,T]\times\bbr^d)$ because $M([0,T]\times\bbr^d)$ is a Banach space. In addition, $\mu_T*\rho_m\to\mu_T*\rho_T$ in $M([0,T]\times\bbr^d)$ by Lemma~\ref{prop1}, so that we get $\rho_T + \mu_T =\mu_T*\rho_T$. 
	
	In the general case where $\|\mu_T\|$ is not necessarily smaller than one, we consider the measure $\la_m(\dd s, \dd y):=\ee^{-ms}\,\mu_T(\dd s, \dd y)$ and note that for sufficiently large $m$ we have $\|\la_m\|<1$ because $\mu(\{0\}\times\bbr^d)=0$. In this case, by what we have already proved, there exists a measure $\eta_m$ satisfying $\eta_m+\la_m=\la_m*\eta_m$. But then $\rho_T(\dd s, \dd y):=\ee^{ms}\,\eta_m(\dd s, \dd y)$ satisfies
	\begin{align*}
	(\rho_T+\mu_T)(\dd s, \dd y)&=\ee^{ms}\,\eta_m(\dd s, \dd y)+\ee^{ms}\ee^{-ms}\,\mu_T(\dd s, \dd y)=\ee^{ms}\,\eta_m(\dd s, \dd y)+\ee^{ms}\,\la_m(\dd s, \dd y)\\
	&=\ee^{ms}\,(\eta_m(\dd s, \dd y)+\la_m(\dd s, \dd y))=\ee^{ms}\,(\la_m*\eta_m)(\dd s, \dd y)\\
	&=([\ee^{ms}\,\la_m(\dd s, \dd y)]*[\ee^{ms}\,\eta_m(\dd s, \dd y)])(\dd s, \dd y)=(\mu_T*\rho_T)(\dd s, \dd y),
	\end{align*}
	where the fifth equation follows from the definition of the convolution. Thus $\rho_T + \mu_T =\mu_T*\rho_T$.
	
	In order to show the uniqueness of $\rho_T$, we assume that there are $\rho_T$ and $\eta_T$ in $M([0,T]\times\bbr^d)$ with $\rho_T + \mu_T =\mu_T*\rho_T$ and $\eta_T + \mu_T =\mu_T*\eta_T$.
	Then
	\begin{align*}	\rho_T=&~\mu_T*\rho_t-\mu_T=(\mu_T*\eta_T-\eta_T)*\rho_T-\mu_T=\eta_T*(\mu_T*\rho_T-\rho_T)-\mu_T=\eta_T*\mu_T-\mu_T\\
	=&~\eta_T. \end{align*}
	
	Now, having constructed $\rho_T$ for every positive $T$ and noting that for every $j\in\bbn$ the restriction of $\rho_{j+1}$ to $[0,j]\times\bbr^d$ must be equal to $\rho_j$ by uniqueness, we define $\rho$ to be the unique measure on $\bbr_+\times\bbr^d$ with $\rho=\rho_T$ on $[0,T]\times\bbr^d$. We still have $\rho \in M_{\loc}(\bbr_+\times\bbr^d)$ and $\rho+\mu=\mu*\rho$, so the proof of (1) is complete. 

	\vspace{0.5\baselineskip}\noindent (2)	Let $\rho$ be the resolvent of $\mu$ as in part (1). Then for $F\in L_{\loc}^p(\bbr_+\times\bbr^d)$ define $X$ by \eqref{eqndetsolution}, which is well defined by Lemma~\ref{prop2}. Also by Lemma~\ref{prop2}, we obtain $X\in L_{\loc}^p(\bbr_+\times\bbr^d)$ and 
	\[ X-\mu*X=X-\mu*(F-\rho*F)=X-(\mu-\mu*\rho)*F=X+\rho*F=F, \]
	thus $X$ is a solution of \eqref{eqnvolterra}.
	To show uniqueness let $\tilde{X}$ be an arbitrary solution of \eqref{eqnvolterra} in $L_{\loc}^p(\bbr_+\times\bbr^d)$. Then
	\[ \tilde{X}=F+\mu*\tilde{X}=F+(\rho*\mu-\rho)*\tilde{X}=F-\rho*(\tilde{X}-\mu*\tilde{X})=F-\rho*F, \]
	hence $\tilde{X}=X$.
	
	\vspace{0.5\baselineskip}\noindent(3) The assumptions on $F$ guarantee that $X$ belongs to $\mathcal{L}(\mu)$ and that all calculations in the previous part remain valid.
	\halmos
\end{Proof}

\section{Examples of resolvents}\label{AppB}

In this section we derive formulae for temporo-spatial resolvents in various examples. They are based on the following lemma whose proof follows directly from the definition of the resolvent measure.
\blem\label{resLemma} Suppose that $\mu\in M_\loc(\bbr_+\times\bbr^d)$ has a resolvent measure $\rho\in M_\loc(\bbr_+\times\bbr^d)$.
\benu
\item If $\mu=m\otimes \delta_{0,\bbr^d}$ with some $m\in M_\loc(\bbr_+)$ satisfying $m(\{0\})=0$ and $r$ is the temporal resolvent measure of $m$ (i.e., the unique $r\in M_\loc(\bbr_+)$ with $r*m=r+m$), then $\rho=r\otimes\delta_{0,\bbr^d}$.
\item If $\mu$ has a Lebesgue density $k\in L^1_\loc(\bbr_+\times\bbr^d)$, then also $\rho$ has a Lebesgue density $r\in L^1_\loc(\bbr_+\times\bbr^d)$, which is given by
\beq\label{r-formula} r(t,x) = -\sum_{n=1}^\infty k^{*n}(t,x), \eeq
where the series converges absolutely for almost every $(t,x)\in\bbr_+\times\bbr^d$.
\item If in the situation of (2) we have that $k(t,x)=-\la f(x)$ with some $f\in L^1(\bbr^d)$ $\la\in\bbr$, then \eqref{r-formula} takes the form
\begin{equation*}
	%\label{r-density} 
	r(t,x) = \la \sum_{n=1}^\infty \frac{(-\la t)^{n-1}}{(n-1)!} f^{*n}(x),\quad (t,x)\in\bbr_+\times\bbr^d. 	
\end{equation*}
\eenu
\elem

\begin{Example}\label{Exres} We present some applications of Lemma~\ref{resLemma}.
\benu
	\item If $\mu=-\la \Leb_{\bbr_+}\otimes \delta_{0,\bbr^d}$ is the measure considered in Examples~\ref{Ex1} and \ref{Ex2}, the resolvent measure is given by $\rho(\dd t,\dd x) = \la\ee^{-\la t}\,\dd t \,\delta_{0,\bbr^d}(\dd x)$.
	%\item Let $f(x)=(2\pi)^{-d/2}\ee^{-|x|^2/2}$ be the density of the $d$-dimensional standard normal distribution. Then $f^{*n}(x)=(2\pi n)^{-d/2}\ee^{-|x|^2/(2n)}$ and
		%\[ r(t,x)=-\sum_{n=1}^\infty\frac{(\pm t)^{n-1}}{(2\pi n)^{d/2}(n-1)!}\ee^{-\frac{|x|^2}{2n}},\quad (t,x)\in\bbr_+\times\bbr. \]
		%There is no analytical expression for this formula. However, this series converges very fast, allowing for efficient numerical computations, see \CC{Fig}.	
		\item Let $f(x)=1/(\pi(1+x^2))$ be the density of the one-dimensional $\mathrm{Cauchy}(0,1)$-distribution. Then $f^{*n}(x)=n/(\pi(x^2+n^2))$ is the density of the $\mathrm{Cauchy}(0,n)$-distribution and
		\[ r(t,x)=\la \sum_{n=1}^\infty\frac{(-\la t)^{n-1}}{(n-1)!}\frac{n}{\pi(x^2+n^2)} = \frac{\la}{2\pi x}G(\la t,x),\quad (t,x)\in\bbr_+\times\bbr, \]
		where $G$ is the (real-valued) function given by
		\begin{align} G(t,x):=&~\ii \Big( t^{-(1+\ii x)}\big(\Ga(1+\ii x)-\Ga(1+\ii x,t)-\Ga(2+\ii x)+\Ga(2+\ii x,t)\big) \nonumber\\
		&~- t^{-(1-\ii x)}\big(\Ga(1-\ii x)-\Ga(1-\ii x,t)-\Ga(2-\ii x)+\Ga(2-\ii x,t)\big)   \Big) \label{funcG} \end{align}
		and $\Ga(\cdot,\cdot)$ is the upper incomplete gamma function.
		\item If $f(x)=\ee^{-x}\bone_{\bbr_+}(x)$, then $f^{*n}(x)=\ee^{-x} x^{n-1}((n-1)!)^{-1}\bone_{\bbr_+}(x)$ and hence
		\[ r(t,x)=\la \sum_{n=1}^\infty\frac{(-\la t)^{n-1}}{(n-1)!}\frac{x^{n-1}\ee^{-x}}{(n-1)!}\bone_{\bbr_+}(x)=\begin{cases} \la \ee^{-x} J_0(\sqrt{2\la t x})\bone_{\bbr_+}(x) &\text{if }\la \geq 0,\\ \la \ee^{-x} I_0(\sqrt{|2\la t x|})\bone_{\bbr_+}(x) &\text{if } \la < 0,  \end{cases}\]
		where $J_0$ ($I_0$) is the (modified) zeroth order Bessel function of the first kind.
		\item If $k$ is a multiple of the heat kernel, that is, $k(t,x)=\la(4\pi t)^{-d/2}\exp(-|x|^2/(4t))\bone_{(0,\infty)}(t)$ for some $\la\in\bbr$, we have $k^{*n}(t,x) = (\la t)^{n-1}/(n-1)!k(t,x)$ and therefore
		\[ r(t,x)=-\sum_{n=1}^\infty \frac{(\la t)^{n-1}}{(n-1)!} k(t,x) = -\ee^{\la t}k(t,x),\quad (t,x)\in\bbr_+\times\bbr^d.\] \halmos
\eenu
	
\end{Example}

\section{Integrability properties of resolvents}\label{AppC}
In many cases integrability properties of resolvents are of interest, see for example Theorem~\ref{stationarysol}. In this section we present criteria for the resolvent $\rho$ to lie in $M(\bbr_+\times\bbr^d)$. The conditions of the lemmata below are given in terms of the Laplace transform of the drift measure $\mu$ which is defined as follows: for a measure $\mu$ in $M(\bbr^d)$ the Laplace transform $\hat{\mu}(z)$ is the function
\[ \hat{\mu}(z)=\int_{\bbr^d}\ee^{-z\cdot u}\, \mu(\dd u), \]
defined for those $z\in\bbc^d$ for which the integral exists and where $z\cdot u$ denotes the standard scalar product in $\bbc^d$. Every $\mu\in M(\bbr^d)$ can be split into three parts, namely the absolutely continuous part $\mu_c$, the discrete part $\mu_d$ and the singular continuous part $\mu_s$. The next lemma can be proved analogously to the one-parameter case (see Theorem~4.4.3 of \cite{Gripenberg90}). We use the notation $\rpart(z)$ for the real part of a complex number $z$ and $\ii\bbr^d$ for the subspace of $\bbc^d$ consisting of all vectors whose entries have all real part zero.
\begin{Lemma}
	Let $\mu\in M(\bbr_+\times\bbr^d)$ satisfy
	\begin{itemize}
		\item $\mu(\{0\}\times \bbr^d)=0$,
		\item $\hat{\mu}(\tau,\xi)\neq1$ for all $\rpart(\tau)\geq0$, $\xi\in \ii\bbr^d$ and
		\item $\inf_{\rpart(\tau)\geq0,\xi\in \ii\bbr^d}|\widehat{\mu_d}(\tau,\xi)-1|>\|\mu_s\|$.
	\end{itemize}
	Then the resolvent $\rho$ of $\mu$ belongs to $M(\bbr_+\times\bbr^d)$.	
\end{Lemma}
In the special case where $\mu$ is absolutely continuous we can give a condition which is both sufficient and necessary. Once again, the proof is analogous to the one-parameter case, cf. Theorem~2.4.1 of \cite{Gripenberg90}.
\begin{Lemma}
	Let $\mu\in M(\bbr_+\times\bbr^d)$ be absolutely continuous. Then the resolvent $\rho$ of $\mu$ belongs to $M(\bbr_+\times\bbr^d)$ if and only if $\hat{\mu}(\tau,\xi)\neq1$ for all $\rpart(\tau)\geq0$ and $\xi\in \ii\bbr^d$.
\end{Lemma}
Integrability of $\rho$ is sufficient, but not necessary for the convolution operator $g\mapsto\rho*g$ to map $L^2_{\loc}(\bbr_+\times\bbr^d)$ into $L^2_{\loc}(\bbr_+\times\bbr^d)$. The following lemma is analogous to the one-parameter case in Theorem~2.6.2 of \cite{Gripenberg90} and provides a criterion in this respect.
\begin{Lemma}
	Let $\mu\in M(\bbr_+\times\bbr^d)$ be absolutely continuous and satisfy the conditions
	\begin{itemize}
		\item $\sup_{\si>0, (\tau,\xi)\in\bbr^d}\left| \frac{\hat{\mu}(\si+\ii\tau,\ii\xi)}{\hat{\mu}(\si+\ii\tau,\ii\xi)-1} \right|<\infty $,
		\item $\int_{\bbr_+\times\bbr^d}\ee^{-\si t}\,|\mu|(\dd t,\dd x)<\infty$ for all $\si>0$.
	\end{itemize}
	Then the convolution operator $g\mapsto\rho*g$ maps $L^2_{\loc}(\bbr_+\times\bbr^d)$ continuously into $L^2_{\loc}(\bbr_+\times\bbr^d)$, where $\rho$ is the resolvent of $\mu$.
\end{Lemma}

\end{appendix}

\subsection*{Acknowledgement}
The authors cordially thank Gustaf Gripenberg, Jean Jacod and Claudia Klüppelberg for inspiring discussions and valuable advice. Furthermore, the first author gratefully acknowledges support from the graduate program TopMath at the Technical University of Munich.

%\addcontentsline{toc}{section}{References}
%\bibliographystyle{plainnat}
%\bibliography{bib-OrnsteinUhlenbeck}

\end{document}